\documentclass[12pt]{amsart}

\usepackage[margin=1.1in]{geometry}
\usepackage[english]{babel}
\usepackage[utf8]{inputenc}
\usepackage{fancyhdr}
\usepackage{natbib}
\usepackage{amsmath}

\RequirePackage{amsthm,amsmath,amsfonts,amssymb}
\usepackage{graphicx}
\usepackage{hyperref}
\usepackage{graphicx}
\usepackage{placeins}
\usepackage{caption}
\usepackage{subcaption}
\numberwithin{equation}{section}

\theoremstyle{plain}
\theoremstyle{remark} 
 \numberwithin{equation}{section} 

\title{Componentwise Equivariant Estimation of Order Restricted Location and Scale Parameters In Bivariate Models: A Unified Study}
\author{Naresh Garg  and Neeraj Misra \\ {\footnotesize Department of Mathematics and Statistics\\Indian Institute of Technology Kanpur \\Kanpur-208016, Uttar Pradesh, India}}

\makeatletter
\def\@seccntformat#1{%
  \protect\textup{\protect\@secnumfont
    \ifnum\pdfstrcmp{subsection}{#1}=0 \bfseries\fi% subsection # in \bfseries
    \csname the#1\endcsname
    \protect\@secnumpunct
  }%
}  
\makeatother

\begin{document}
\maketitle
\section*{\textbf{Abstract}}

	The problem of estimating location (scale) parameters $\theta_1$ and $\theta_2$ of two distributions when the ordering between them is known apriori (say, $\theta_1\leq \theta_2$) has been extensively studied in the literature. Many of these studies are centered around deriving estimators that dominate the best location (scale) equivariant estimators, for the unrestricted case, by exploiting the prior information that $\theta_1 \leq \theta_2$. Several of these studies consider specific distributions such that the associated random variables are statistically independent. This paper considers a general bivariate model and general loss function and unifies various results proved in the literature. We also consider applications of these results to a bivariate normal and a Cheriyan and Ramabhadran's bivariate gamma model. A simulation study is also considered to compare the risk performances of various estimators under bivariate normal and Cheriyan and Ramabhadran's bivariate gamma models.
\\~\\ \textbf{Keywords:} Best location equivariant estimator (BLEE); Best scale equivariant estimator (BSEE), Brewster-Zidek type estimator; Generalized Bayes estimators; Stein type estimator.

\section{\textbf{Introduction}}\label{intro}
\setcounter{equation}{0}
\renewcommand{\theequation}{1.\arabic{equation}}
Let $\underline{X}=(X_1,X_2)$ be a random vector having a joint probability density function (pdf) belonging to location (scale) family
\begin{equation}\label{eq:1.1}
f_{\underline{\theta}}(x_1,x_2)= 	f(x_1-\theta_1,x_2-\theta_2),\; \; \;(x_1,x_2)\in \Re^2, 
\end{equation} 
\begin{equation}\label{eq:1.2}
\left(\quad	f_{\underline{\theta}}(x_1,x_2)= 	\frac{1}{\theta_1 \theta_2}f\!\left(\frac{x_1}{\theta_1},\frac{x_2}{\theta_2}\right),\;\;\; \; \;(x_1,x_2)\in \Re^2 \quad \right),\;
\end{equation} 
where $\underline{\theta}=(\theta_1,\theta_2)\in \Theta=\Re^2 \,(\Re_{++}^2)$ is the vector of unknown location (scale) parameters and $f(\cdot,\cdot) $ is a specified pdf on $\Re^2$; here $\Re^2=(-\infty,\infty)\times(-\infty,\infty)$ and $\Re_{++}^2=(0,\infty)\times(0,\infty)$. Generally, $\underline{X}=(X_1,X_2)$ would be a minimal-sufficient statistic based on a bivariate random sample or two independent random samples, as the case may be. In many real life situations ordering between the parameters $\theta_1$ and $\theta_2$ may be known apriori (say, $\theta_1\leq\theta_2$) and it may be of interest to estimate $\theta_1$ and $\theta_2$ (see, for example, Barlow et al. (\citeyear{MR0326887}), Robertson et al. (\citeyear{MR961262}), Kumar and Sharma (\citeyear{MR981031}), Kubokawa and Saleh (\citeyear{MR1370413}), Hwang and Peddada (\citeyear{MR1272076}) and references cited therein).\vspace{1mm}

Let $\Theta_0=\{\underline{\theta}\in\Theta: \theta_1\leq\theta_2\} $ be the restricted parameter space. There is an extensive literature on estimation of $\theta_1$ and $\theta_2$ (simultaneously, as well as, componentwise) when it is known apriori that $\underline{\theta}\in \Theta_0$. A natural question that arises in these problems is whether the best location (scale) equivariant estimator(s) (BLEE (BSEE)) can be improved by exploiting the prior information that $ \underline{\theta}\in\Theta_0 $. Many researchers have studied this and related aspects of the problem. However, several of these studies are focussed to specific distributions, having independent marginals, and specific loss functions. Some of the contributions in this direction are due to Cohen and Sackrowitz (\citeyear{MR270483}), Brewster and Zidek (\citeyear{MR381098}), Lee (\citeyear{MR615447}), Kumar and Sharma (\citeyear{MR981031}, \citeyear{MR1058934}), Kelly (\citeyear{MR994278}), Kushary and Cohen (\citeyear{MR1029476}), Gupta and Singh (\citeyear{MR1210354}), Pal and Kushary (\citeyear{MR1165709}), Misra and Singh (\citeyear{MR1366828}), Vijayasree et al. (\citeyear{MR1345425}), Misra and Dhariyal (\citeyear{MR1326266}), Misra et al. (\citeyear{MR1904424}, \citeyear{MR2205815}), Oono and Shinozaki (\citeyear{OonoY}), Chang and Shinozaki (\citeyear{MR3390170}), Petropoulos (\citeyear{Petropoulos}, \citeyear{PetropoulosC}), Bobotas (\citeyear{MR3963133}, \citeyear{MR3926974}) and Patra et al. (\citeyear{doi:10.1080/02331888.2021.1943395}). For a few contributions to this problem under general setting (general probability model and/or general loss function) readers may refer to  Blumenthal and Cohen (\citeyear{MR223007}), Sackrowitz (\citeyear{MR254961}), Hwang and Peddada (\citeyear{MR1272076}), Kubokawa and Saleh (\citeyear{MR1370413}) and Iliopoulos (\citeyear{MR1804623}). For a detailed account of contributions in this area of research one may refer to the research monograph by van Eeden (\citeyear{MR2265239}).
\\~\\ Kubokawa and Saleh (\citeyear{MR1370413}) considered the location (scale) model (1.1) ((1.2)) with
$$f(z_1,z_2)=f_1(z_1)f_2(z_2),\;\; -\infty<z_i<\infty,\;i=1,2,$$
where $f_1$ and $f_2$ are specified pdfs on the real line $\Re$. They dealt with estimation of the smaller location (scale) parameter $\theta_1$ when it is known apriori that  $\underline{\theta}\in\Theta_0$. They assumed that 
\begin{small}
\begin{equation}
	f_i\in P_L=\bigg\{g\,:\, \frac{g(y+c_1)}{g(y+c_2)} \text{ is non-decreasing in } y\in \Re,\text{ for every }c_1<c_2\bigg\},\;i=1,2,
\end{equation}
\begin{equation} \label{eqn:(1.4)}
	\left(	f_i\in P_S=\bigg\{g\,:\,\frac{g(c_1y)}{g(c_2y)} \text{ is non-decreasing in } y\in \Re_{++},\text{ for every }0<c_1<c_2\bigg\},\;i=1,2, \right)
\end{equation}
\end{small}
and considered a quite general loss function 
$$L_1(\underline{\theta},a)=W(a-\theta_1),  \;\underline{\theta}\in\Theta_0,\; a\in\Re,\; \left(L_1(\underline{\theta},a)=W\left(\frac{a}{\theta_1}\right),  \;\underline{\theta}\in\Theta_0,\; a\in\Re_{++}\right),$$
where $W:\Re\rightarrow[0,\infty)$ is such that $W(0)=0$ ($W(1)=0$), $W(t)$ is strictly decreasing for $t<0$ ($t<1$)  and strictly increasing for $t\geq0$ ($t\geq 1$). Under the above set-up, they derived conditions that ensure improvements over the best location (scale) equivariant estimator of $\theta_1$. They found explicit expressions of the dominating estimators. In fact, Kubokawa and Saleh (\citeyear{MR1370413}) dealt with estimation of the smallest location (scale) parameter of $k\;(\,\geq2)$ independent location (scale) families of probability distributions when it is known apriori that the corresponding location (scale) parameters $\theta_1,\theta_2,\!...,\theta_k$ satisfy the tree ordering ($\theta_1\leq\theta_i,\; i=2,3,...,k$). Under the set-up considered by Kubokawa and Saleh (\citeyear{MR1370413}), the random variables $X_1$ and $X_2$ are independently distributed. In this paper, we extend the study of Kubokawa and Saleh (\citeyear{MR1370413}) to situations where $X_1$ and $X_2$ may be statistically dependent. As in Kubokawa and Saleh (\citeyear{MR1370413}), we will closely follow the IERD (Integral expression risk difference) approach of Kubokawa (\citeyear{MR1272084}), to obtain improvements over the BLEE/BSEE of $\theta_1$ and $\theta_2$. We also consider estimation of the larger location (scale) parameter $\theta_2$ that has not been addressed by Kubokawa and Saleh (\citeyear{MR1370413}). To avoid some notational and presentation difficulties, throughout the paper, we extend the usual orders "$\leq$" and "$<$" in the real line to the extended real line $\Re\cup \{-\infty,\infty\}$ with the following convention. For any positive (negative) real number "$b$", we take $\frac{b}{0}=\infty\, (-\infty)$ and, for any real number "$c$", we take $-\infty<c<\infty$.   \vspace{1mm}

In Section \ref{sec:2} (\ref{sec:3}), we consider a general bivariate location (scale) family of distributions and deal with componentwise estimation of order restricted location (scale) parameters $\theta_1$ and $\theta_2$ under a quite general loss function. We derive sufficient conditions that guarantee improvements over the BLEE (BSEE). The explicit expressions of dominating estimators are obtained. In Subsection \ref{sec:2.3} (\ref{sec:3.3}), We provide applications of various results derived in the paper. In Subsection \ref{sec:2.4} (\ref{sec:3.4}), we consider a simulation study for comparing risk performances of various estimators of smaller location (scale) parameter under bivariate normal (Cheriyan and Ramabhadran's bivariate gamma) model.

\section{\textbf{Improving the Best Location Equivariant Estimators (BLEEs)}}
\label{sec:2}
\setcounter{equation}{0}
\renewcommand{\theequation}{2.\arabic{equation}}

Firstly, we will introduce some notations in connection with the probability model \eqref{eq:1.1}. Let $Z_i=X_i-\theta_i$, $i=1,2,$ and $\underline{Z}=(Z_1,Z_2),$ so that $\underline{Z}$ has the joint pdf $f(t_1,t_2), \; (t_1,t_2)\in\Re^2$. Let $S_i$ be the distributional support of $Z_i=X_i-\theta_i,\;i=1,2$. Let $Z=Z_2-Z_1$ and $f_i$ be the pdf of $Z_i,\; i=1,2$. Then, 
\\$f_1(s)=\int\limits_{-\infty}^{\infty} f(s,t) \,dt,\;\; s\in \Re, \qquad f_2(s)=\int\limits_{-\infty}^{\infty} f(t,s) \,dt,\;\; s\in \Re.$\vspace*{2mm}

For any $s\in S_i$, let $Z_s^{(i)}$ denote a random variable having the same distribution as conditional distribution of $Z$ given $Z_i=s,\;i=1,2$. Then the pdf and the distribution function (df) of $Z_s^{(1)}\,(s\in S_1)$ are given by
\begin{equation*}
h_1(t|s)=\frac{f(s,t+s)}{f_1(s)}, \;\; t\in\Re, \text{ and }\; H_1(t|s)=\frac{\int_{-\infty}^{t}f(s,z+s)dz}{f_1(s)}, \;\; t\in\Re,
\end{equation*}
respectively, and the pdf and the df of $Z_s^{(2)}\,(s\in S_2)$ are given by
\begin{equation*}
h_2(t|s)=\frac{f(t-s,s)}{f_2(s)}, \;\; t\in\Re, \text{ and }\; H_2(t|s)=\frac{\int_{-\infty}^{t}f(z-s,s)dz}{f_2(s)}, \;\; t\in\Re,
\end{equation*}
respectively.
\\~\\For the location model \eqref{eq:1.1}, consider estimation of $\theta_i$ under the loss function 
\begin{equation}\label{eq:2.1}
L_i(\underline{\theta},a)=W(a-\theta_i),\; \underline{\theta}\in\Theta,\; \; a\in \mathcal{A}=\Re,\;i=1,2,
\end{equation}  
where $W:\Re\rightarrow [0,\infty)$ is a specified non-negative function. We make the following assumptions on the function $W(\cdot)$:
\\~\\ \textbf{\boldmath$A_1$:} $W:\Re \rightarrow [0,\infty)$ is such that $W(0)=0$, $W(t)$ is decreasing on $(-\infty,0)$ and increasing on $(0,\infty)$. Further $W^{'}(t)$ is non-decreasing on the set $D_0$ (the set of points at which $W(\cdot)$ is differentiable). \vspace{2mm}

First consider estimation of $\theta_i,\,i=1,2,$ under the unrestricted parameter space $\Theta=\Re^2$ and the loss function \eqref{eq:2.1}. Under the above set-up, the problem of estimating $\theta_i\;(i=1,2)$ is invariant under the additive group of transformation $\mathcal{G}=\{g_{c_1,c_2}:(c_1,c_2)\in\Re^2 \}$, where $g_{c_1,c_2}(x_1,x_2)=(x_1+c_1,x_2+c_2)$, $(x_1,x_2)\in\Re^2,\;  (c_1,c_2)\in\Re^2$. Any (non-randomized) location equivariant estimator $\delta_i$ of $\theta_i$ is of the form
$ \delta_{c,i}(X_1,X_2)=X_i-c, \; \; i=1,2,$
for some constant $c\in\Re$. The risk function of $\delta_{c,i}$ is given by
$R_i(\underline{\theta},\delta_{c,i})=E_{\underline{\theta}}[L_i(\underline{\theta},\delta_{c,i}(\underline{X}))],\; \, \underline{\theta}\in\Theta,\; i=1,2.$\vspace*{2mm}

The risk function of any location (scale) invariant estimator $\delta_{c,i}$ of $\theta_i$ is constant (does not depend on $\underline{\theta}\in\Theta$). For the existence of unrestricted ($\underline{\theta}\in\Theta$) BLEE, we need the following assumption:
\\~\\ \textbf{\boldmath$A_2$:} The equation $E[W^{'}(Z_i-c)]=0$ has the unique solution, say $c=c_{0,i},\;i=1,2.$ \vspace{2mm}

Since the risk function of any equivariant estimator of $\theta_i$ is constant on $\Theta$, under assumptions $A_1$ and $A_2$, the unique BLEE of $\theta_i$ is 
\begin{equation}\label{eq:2.2}
\delta_{c_{0,i},i}(\underline{X})=X_i-c_{0,i},\;\; i=1,2,
\end{equation}
where $c_{0,i}$ is the unique solution of the equation 
\begin{equation}\label{eq:2.3}
\int\limits_{-\infty}^{\infty} \, W^{'}(z-c)\,f_i(z)dz=0,\;i=1,2.
\end{equation}

Now consider estimation of location parameter $\theta_i, \;i=1,2,$ under the restricted parameter space $\Theta_0=\{(x_1,x_2)\in \Theta: \,x_1 \leq x_2\}$ and the loss function \eqref{eq:2.1}. Under the restricted parameter space $\Theta_0$, the location family of distributions \eqref{eq:1.1} is not invariant under the group of transformations $\mathcal{G}=\{g_{c_1,c_2}:(c_1,c_2)\in\Re^2 \}$, considered above. An appropriate group of transformations ensuring invariance under restricted parameter space $\Theta_0$ is $\mathcal{K}=\{k_c:\,c\in\Re\},$ where $k_c(x_1,x_2)=(x_1+c,x_2+c),\; (x_1,x_2)\in\Re^2,\;c\in\Re$. Under the group of transformations $\mathcal{K}$, the problem of estimating $\theta_i$, under $\underline{\theta}\in\Theta_0$ and the loss function \eqref{eq:2.1}, is invariant. Any location equivariant estimator of $\theta_i$ is of the form 
\begin{equation}\label{eq:2.4}
\delta_{\psi_i}(\underline{X})=X_i-\psi_i(D),
\end{equation}
for some function $\psi_i:\,\Re\rightarrow \Re\,,\;i=1,2,$ where $D=X_2-X_1$. Here the risk function
\begin{equation}\label{eq:2.5}
R_i(\underline{\theta},\delta_{\psi_i})=E_{\underline{\theta}}[L_i(\underline{\theta},\delta_{\psi_i}(\underline{X}))],\; \, \underline{\theta}\in\Theta_0,
\end{equation}
of any location equivariant estimator $\delta_{\psi_i}$ of $\theta_i,\;i=1,2,$ may not be constant on $\Theta_0$, and it depends on $\underline{\theta}\in\Theta_0$ only through $\lambda=\theta_2-\theta_1\in[0,\infty)$. \vspace*{2mm}

The following lemma will be useful in proving the main results of the paper. The proof of the lemma is straight forward and hence omitted. 
\\~\\ \textbf{Lemma 2.1.} Let $s_0\in \Re$ and let $M:\Re\rightarrow\Re$ be such that $M(s)\leq 0,\; \forall \; s<s_0, $ and $M(s)\geq 0,\; \forall \; s> s_0$. Let $ M_i:\Re\rightarrow [0,\infty), \; i=1,2,$ be non-negative functions such that
$M_1(s) M_2(s_0) \geq (\leq)\, M_1(s_0) M_2(s),\; \forall \; s<s_0,
\text{ and } M_1(s) M_2(s_0) \leq\,(\geq)\; M_1(s_0) M_2(s),\; \forall \; s$ $>s_0.$
Then, 
$$ M_2(s_0) \int\limits_{-\infty}^{\infty} M(s) \, M_1(s) ds\leq\;(\geq)\; M_1(s_0) \int\limits_{-\infty}^{\infty} M(s) \, M_2(s) ds.$$

The facts stated in the following lemma are well known in the theory of stochastic orders (see Shaked and Shanthikumar (\citeyear{MR2265633})). The proof of the lemma is straight forward, hence skipped.
\\~\\ \textbf{Lemma 2.2.} If, for any fixed $\Delta\geq 0$ and $t\in \Re$, $h_i(t-\Delta\vert s)/h_i(t\vert s)$ is non-decreasing (non-increasing) in $s\in S_i$,
then $H_i(t-\Delta\vert s)/H_i(t\vert s)$ is non-decreasing (non-increasing) in $s\in S_i$ and $h_i(t\vert s)/H_i(t\vert s)$ is non-increasing (non-decreasing) in $s\in S_i,\;i=1,2.$

\vspace*{2mm}

In the next subsection, we consider equivariant estimation of location parameter $\theta_1$ under the loss function $L_1$, defined by (2.1), when it is known apriori that $\underline{\theta}\in\Theta_0$. We aim to find estimators that dominate the BLEE $\delta_{c_{0,i},i}(\underline{X}),\,i=1,2$ (defined through \eqref{eq:2.2} and \eqref{eq:2.3}) by exploiting the prior information that $\underline{\theta}\in\Theta_0$.
\subsection{\textbf{Improvements Over the BLEE of $\theta_1$}}
\label{sec:2.1}
\setcounter{equation}{0}
\renewcommand{\theequation}{2.1.\arabic{equation}}

\noindent
\vspace*{2mm}

Consider estimation of $\theta_1$ under the loss function $L_1(\underline{\theta},a)=W(a-\theta_1),\;\underline{\theta}\in\Theta_0,\;a\in \Re$, when it is known apriori that $\underline{\theta}\in\Theta_0$. Throughout this subsection, we will assume that the function $W(\cdot)$ satisfies assumptions $A_1$ and $A_2$. 
\vspace*{2mm}

In the following theorem, we provide a class of estimators that improve upon the BLEE $\delta_{c_{0,1},1}(\underline{X})=X_1-c_{0,1}$, defined by \eqref{eq:2.2} and \eqref{eq:2.3}.
\\~\\ \textbf{Theorem 2.1.1.} Suppose that, for any fixed $\Delta\geq 0$ and $t$, $H_1(t-\Delta\vert s)/H_1(t\vert s)$ is non-decreasing (non-increasing) in $s\in S_1$. Consider a location equivariant estimator $\delta_{\psi_1}(\underline{X})=X_1 - \psi_1(D)$ of $\theta_1$, where $\psi_1(t)$ is non-increasing (non-decreasing) in $t$, $\lim_{t\to\infty} \psi_1(t)=c_{0,1}$ and $\int_{-\infty}^{\infty} W^{'}(s-\psi_1(t))\; H_1(t\vert s)\,f_1(s)ds \, \geq \,(\leq) \,0,\; \forall\; t.$ Then $$R_1(\underline{\theta},\delta_{\psi_1})\leq R_1(\underline{\theta},\delta_{c_{0,1},1})\, \;\;\; \forall \; \; \underline{\theta} \in \Theta_0.$$
\begin{proof}
Let us fix $\underline{\theta}\in\Theta_0$ and let $\lambda=\theta_2-\theta_1$, so that $\lambda \geq 0$. Consider the risk difference
\begin{align*}
	\Delta_1(\lambda)&= R_1(\underline{\theta},\delta_{c_{0,1},1})-R_1(\underline{\theta},\delta_{\psi_1})\\
	& = E_{\underline{\theta}}[W(Z_1-c_{0,1})- W(Z_1-\psi_1(Z+\lambda))] \\
	&= E_{\underline{\theta}}\left[\int_{Z+\lambda}^{\infty}\Big\{ \frac{d}{dt} W(Z_1-\psi_1(t))\Big\}\; dt\right]\\
	&= -\int_{-\infty}^{\infty} \psi_1^{'}(t) E_{\underline{\theta}}[W^{'}(Z_1-\psi_1(t))\; I_{(-\infty,t-\lambda]}(Z)\;]\,dt,
\end{align*}
where, for any set $A$, $I_A(\cdot)$ denotes its indicator function. Since $\psi_1(t)$ is a non-increasing (non-decreasing) function of $t$, it suffices to show that, for every $t$ and $\lambda\geq 0$,
\begin{small}
	\begin{align}\label{eq:2.1.1}
		E_{\underline{\theta}}[W^{'}(Z_1-\psi_1(t))\; I_{(-\infty,t-\lambda]}(Z)\;]
		&=E_{\underline{\theta}}[W^{'}(Z_1-\psi_1(t))\, H_1(t-\lambda\vert Z_1)]\nonumber\\
		&=\int_{-\infty}^{\infty} W^{'}(s-\psi_1(t))\;H_1(t-\lambda\vert s)f_1(s)ds\; \geq\;(\leq) \; 0.
	\end{align}
\end{small}
To prove \eqref{eq:2.1.1}, let us fix $t$ and $\lambda\geq 0$. Let $s_0 =\psi_1(t),\; M_1(s)=H_1(t\vert s),\; M_2(s)=H_1(t-\lambda\vert s)$ and $M(s)=W^{'}(s-s_0)f_1(s),\; s\in \Re$. Then, using $A_1,$ we have $M(s)\leq 0,\; \forall \; s<s_0$ and $M(s)\geq 0,\; \forall \; s> s_0.$ Also, under the hypothesis of the theorem,
\begin{align*}
	M_1(s) M_2(s_0) &\geq\;(\leq) \; M_1(s_0) M_2(s),\; \forall \; s<s_0\\
	\text{and }\; \; M_1(s) M_2(s_0) &\leq\;(\geq)\; M_1(s_0) M_2(s),\; \forall \; s>s_0.
\end{align*}
Using Lemma 2.1, we get
\begin{small}
	\begin{align*}
		0\leq\;(\geq)\;H_1(t-\lambda\vert \psi_1(t)) \int_{-\infty}^{\infty}\; W^{'}(s-\psi_1(t))\; H_1(t\vert s)\,f_1(s)  ds \qquad \qquad \qquad \qquad \qquad \qquad \qquad \quad \\ \qquad \qquad \qquad \qquad \qquad \qquad \leq\;(\geq)\; H_1(t\vert  \psi_1(t)) \int_{-\infty}^{\infty}\; W^{'}(s-\psi_1(t))\; H_1(t-\lambda\vert  s)\,f_1(s) ds,
	\end{align*}
\end{small}
which, in turn, implies (2.1.1).
\end{proof}

Now we will prove two useful corollaries to the above theorem. The following corollary provides the Brewster-Zidek (1974) type (B-Z type) improvement over the BLEE $\delta_{c_{0,1},1}$.
\\~\\ \textbf{Corollary 2.1.1. (i)} Suppose that assumptions of Theorem 2.1.1 hold. Further suppose that, for every fixed $t$, the equation
$$k_1(c\vert  t)=\int_{-\infty}^{\infty} \; W^{'}(s-c)\; H_1(t\vert  s)\,f_1(s)\;ds =0$$
has the unique solution $c\equiv \psi_{0,1}(t)\in S_1$. Then
$$R_1(\underline{\theta},\delta_{\psi_{0,1}})\leq R_1(\underline{\theta},\delta_{c_{0,1},1})\, \;\;\; \forall \; \; \underline{\theta} \in \Theta_0,$$
where $\delta_{\psi_{0,1}}(\underline{X})=X_1-\psi_{0,1}(D)$.
\\~\\\textbf{(ii)} In addition to assumptions of (i) above, suppose that $\psi_{1,1}:\Re \rightarrow \Re $ is such that $ \psi_{1,1}(t) \leq \; (\geq) \; \psi_{0,1}(t), \; \; \forall \; t, \; \psi_{1,1}(t) $ is non-increasing (non-decreasing) in $t$ and $\lim_{t \to \infty}\; \psi_{1,1}(t) = c_{0,1}$.
Then 
$$R_1(\underline{\theta}, \delta_{\psi_{1,1}}) \leq R_1(\underline{\theta}, \delta_{c_{0,1},1}),\; \; \forall \; \underline{\theta} \in \Theta_0,$$
where $\delta_{\psi_{1,1}}(\underline{X})=X_1-\psi_{1,1}(D)$.
\begin{proof}
It suffices to show that $\psi_{0,1}(t)$ satisfies conditions of Theorem 2.1.1. Note that the hypothesis of the corollary, along with the assumption $A_2$, ensure that $\lim_{t\to\infty} \psi_{0,1}(t)=c_{0,1}$. To show that $\psi_{0,1}(t)$ is a non-increasing (non-decreasing) function of $t$, suppose that, there exist numbers $t_1$ and $t_2$ such that $t_1<t_2$ and $\psi_{0,1}(t_1)\neq \psi_{0,1}(t_2).$ We have $k_1(\psi_{0,1}(t_1)\vert t_1)=0$. Also, using the hypotheses of the corollary and the assumption $A_1$, it follows that $\psi_{0,1}(t_2)$ is the unique solution of $k_1(c\vert t_2)=0$ and $k_1(c\vert t_2)$ is a non-increasing function of $ c $. Let $ s_0 = \psi_{0,1}(t_1), \; M(s)=W^{'}(s-s_0)f_1(s),\; M_1(s)=H_1(t_2\vert s)$ and $M_2(s)=H_1(t_1\vert s),\; s\in S_1$. Then, under assumption $A_1$, using Lemma 2.1, we get
\begin{small}
	$$ H_1(t_1\vert \psi_{0,1}(t_1))\; \int_{-\infty}^{\infty}\; W^{'}(s-\psi_{0,1}(t_1))\; H_1(t_2\vert s) f_1(s) ds \qquad \qquad \qquad \qquad \qquad \qquad \quad$$   $$\qquad \qquad \qquad \qquad \qquad \qquad \leq\; (\geq)\; H_1(t_2\vert \psi_{0,1}(t_1))\; \int_{-\infty}^{\infty}\; W^{'}(s-\psi_{0,1}(t_1))\; H_1(t_1\vert s)\,f_1(s)  ds =0.$$
\end{small}
$$\implies \;\; k_1(\psi_{0,1}(t_1)\vert t_2)=\int_{-\infty}^{\infty}\; W^{'}(s-\psi_{0,1}(t_1))H_1(t_2\vert s)\,f_1(s) ds\leq \; (\geq) \; 0.$$
This implies that $k_1(\psi_{0,1}(t_1)\vert t_2)\,< \, (>) \, 0$, as $k_1(c\vert t_2)=0$ has the unique solution $c\equiv \psi_{0,1}(t_2)$ and $\psi_{0,1}(t_1)\neq\psi_{0,1}(t_2)$. Since $k_1(c\vert t_2)$ is a non-increasing function of c, $k_1(\psi_{0,1}(t_2)\vert t_2)$ $=0$ and $k_1(\psi_{0,1}(t_1)\vert t_2)\,< \, (>) \, 0$, it follows that $\psi_{0,1}(t_1)>(<) \psi_{0,1}(t_2)$. 
\\~\\ The proof of part (ii) is an immediate by-product of Theorem 2.1.1 using the fact that, for any $t$, $k_1(c\vert t)$ is a non-increasing function of $c\in\Re$.
\end{proof}

In the following corollary we provide the Stein (\citeyear{MR171344}) type improvements over the BLEE $\delta_{c_{0,1},1}(\underline{X}).$
\\~\\ \textbf{Corollary 2.1.2. (i)} Suppose that, for any fixed $\Delta\geq 0$ and $t$, $h_1(t-\Delta\vert s)/h_1(t\vert s)$ is non-decreasing (non-increasing) in $s\in S_1$. Let $\psi_{0,1}(t)\in S_1$ be as defined in Corollary 2.1.1. In addition suppose that, for any $t$, the equation
$$k_2(c\vert t)=\int_{-\infty}^{\infty} \; W^{'}(s-c)\;h_1(t\vert s)\,f_1(s)\;ds =0$$
has the unique solution $c\equiv \psi_{2,1}(t)\in S_1.$ Let $\psi_{2,1}^{*}(t)= \max\{c_{0,1},\psi_{2,1}(t)\}$ ($\psi_{2,1}^{*}(t)= \min\{c_{0,1},\psi_{2,1}(t)\}$) and $\delta_{\psi_{2,1}^{*}}(\underline{X}) = X_1- \psi_{2,1}^{*}(D)$. Then
$$R_1(\underline{\theta},\delta_{\psi_{2,1}^{*}})\leq R_1(\underline{\theta},\delta_{c_{0,1},1})\, \;\;\; \forall \; \; \underline{\theta} \in \Theta_0.$$
\textbf{(ii)} In addition to assumptions of (i) above, suppose that $\psi_{3,1}:\Re \rightarrow \Re $ be such that $ \psi_{3,1}(t) \leq \, (\geq) \; \psi_{2,1}(t), \; \; \forall \; t$ and $\psi_{3,1}(t) $ is non-increasing (non-decreasing) in $t$. Define $\psi_{3,1}^{*}(t)=\max\{c_{0,1},\psi_{3,1}(t)\}\;(\psi_{3,1}^{*}(t)=\min\{c_{0,1},\psi_{3,1}(t)\})$ and $\delta_{\psi_{3,1}^{*}}(\underline{X})=X_1-\psi_{3,1}^{*}(D)$.
Then 

$$R_1(\underline{\theta}, \delta_{\psi_{3,1}^{*}}) \leq R_1(\underline{\theta}, \delta_{c_{0,1},1}),\; \; \forall \; \underline{\theta} \in \Theta_0.$$
\begin{proof}
It suffices to show that $\psi_{2,1}^{*}(\cdot)$ satisfies conditions of Theorem 2.1.1. Under the assumption that, for any fixed $\Delta\geq 0$ and $t$, $h_1(t-\Delta\vert s)/h_1(t\vert s)$ is non-decreasing (non-increasing) in $s\in S_1$, on following the line of arguments used in proving Corollary 2.1.1, it can be concluded that $\psi_{2,1}(t)$ (and hence $\psi_{2,1}^{*}(t)$) is non-increasing (non-decreasing) in $t$. To show that $\lim_{t \to \infty} \; \psi_{2,1}^{*}(t)=c_{0,1}$, we will show that $\psi_{2,1}(t)\leq\; (\geq)\; \psi_{0,1}(t),\; \forall \; t$. Let us fix $t$. Then 
$ k_1(\psi_{0,1}(t)\vert t)=k_2(\psi_{2,1}(t)\vert t)=0$.
\vspace*{2mm}

The hypothesis of the theorem and Lemma 2.2, imply that, for every fixed $t$, $h_1(t\vert s)/H_1(t\vert s)$ is non-increasing (non-decreasing) in $s\in S_1$. Let $s_0=\psi_{0,1}(t),\; M(s)=W^{'}\!(s-s_0)f_1(s),\; $ $M_1(s)=h_1(t\vert s)$ and $M_2(s)=H_1(t\vert s),\; s\in S_1.$ Using assumption $A_1$, the monotonicity of  $h_1(t\vert s)/H_1(t\vert s)$, Lemma 2.1 and the fact that $k_1(\psi_{0,1}(t)\vert t)=0$, we conclude that
\begin{multline*}
	H_1(t\vert\psi_{0,1}(t))\; \int_{-\infty}^{\infty}\; W^{'}(s-\psi_{0,1}(t))\; h_1(t\vert s) f_1(s)  ds \\ \leq\; (\geq)\;h_1(t\vert \psi_{0,1}(t))\; \int_{-\infty}^{\infty}\; W^{'}(s-\psi_{0,1}(t))\; H_1(t\vert s)\,f_1(s)ds =0
\end{multline*}
\begin{equation}\label{eq:2.1.3}
	\implies \qquad	k_2(\psi_{0,1}(t)\vert t)=\int_{-\infty}^{\infty}\; W^{'}(s-\psi_{0,1}(t))\,h_1(t\vert s)\,f_1(s) ds\leq \; (\geq) \; 0,
\end{equation}
Since $k_2(c\vert t)$ is a non-increasing function of c and $\psi_{2,1}(t)$ is the unique solution of $k_2(c\vert t)=0$, using \eqref{eq:2.1.3}, we conclude that $\psi_{0,1}(t)\geq\,(\leq)\psi_{2,1}(t)$. Hence $c_{0,1}=\lim_{t \to \infty} \psi_{0,1}(t)\geq\, (\leq) \\\; \lim_{t \to \infty}\psi_{2,1}(t)$ and $\lim_{t \to \infty} \psi_{2,1}^{*}(t)=\max\{c_{0,1},\lim_{t \to \infty} \psi_{2,1}(t)\}=c_{0,1}$  ($\lim_{t \to \infty} \psi_{2,1}^{*}(t)= \min\{c_{0,1},\lim_{t \to \infty}\psi_{2,1}(t)\}=c_{0,1}$).  Note that $\psi_{2,1}^{*} (t)\leq (\geq) \psi_{0,1}(t),\; \forall\; t$. Since $k_1(c\vert t)$ is a non-increasing function of $c$, we have 
$$k_1(\psi_{2,1}^{*}(t)\vert t)\geq\;(\leq)\; k_1(\psi_{0,1}(t)\vert t)=0,\; \forall\; t.$$
Hence the result follows. 
\end{proof}

The proof of part (ii) of Corollary 2.1.2 is immediate from Theorem 2.1.1 on noting that $\psi_{0,1}(t)\,\geq \,(\leq)\, \psi_{2,1}(t),\;\forall \; t,$ and $k_1(c\vert t)$ is a non-increasing function of $c$, for every $t$.
\\~\\ \textbf{Remark 2.1.1.}	It is straightforward to see that the Brewster-Zidek (1974) type estimator $\delta_{\psi_{0,1}}$, derived in Corollary 2.1.1 (i), is the generalized Bayes estimator with respect to the non-informative prior density $\pi(\theta_1,\theta_2)=1,\;(\theta_1,\theta_2)\in\Theta_0.$\vspace*{2mm}

The results reported in Theorem 2.1.1, Corollary 2.1 (i) and Corollary 2.2 (i) are extensions of results proved by Kubokawa and Saleh (\citeyear{MR1370413}) for the special case when $X_1$ and $X_2$ are independently distributed. The results for estimating the larger location parameter $\theta_2$ can be obtained along the same lines. For brevity, in the following section, we state these results without providing their proofs.

\subsection{\textbf{Improvements Over the BLEE of $\theta_2$}}
\label{sec:2.2}
\setcounter{equation}{0}
\renewcommand{\theequation}{2.2.\arabic{equation}}

\noindent
\vspace*{2mm}

Under assumptions $A_1$ and $A_2$, consider estimation of $\theta_2$ under the loss function $L_2(\underline{\theta},a)=W(a-\theta_2),\;\underline{\theta}\in\Theta_0,\;a\in \Re$, when it is known apriori that $\underline{\theta}\in\Theta_0$.
\vspace*{2mm}

The following theorem provides a class of estimators that improve upon the BLEE, $\delta_{c_{0,2},2}(\underline{X})=X_2-c_{0,2}$, of $\theta_2$, defined by \eqref{eq:2.2} and \eqref{eq:2.3}.
\\~\\ \textbf{Theorem 2.2.1.} Suppose that, for any fixed $\Delta\geq 0$ and $t$, $H_2(t-\Delta\vert s)/H_2(t\vert s)$ is non-increasing (non-decreasing) in $s\in S_2$. Consider a location equivariant estimator $\delta_{\psi_2}(\underline{X})=X_2 - \psi_2(D)$ of $\theta_2$, where $\psi_2(t)$ is non-decreasing (non-increasing) in $t$, $\lim_{t\to\infty} \psi_2(t)=c_{0,2}$ and $\int_{-\infty}^{\infty} W^{'}(s-\psi_2(t))\; H_2(t\vert s)\,f_2(s)ds \, \leq \,(\geq) \,0,\; \forall\; t.$ Then $$R_2(\underline{\theta},\delta_{\psi_2})\leq R_2(\underline{\theta},\delta_{c_{0,2},2}),\, \;\;\; \forall \; \; \underline{\theta} \in \Theta_0.$$

The following corollary provides the B-Z type improvements over the BLEE $\delta_{c_{0,2},2}(\cdot)$.
\\~\\ \textbf{Corollary 2.2.1.} Suppose that assumptions of Theorem 2.2.1 hold. Further suppose that, for every fixed $t$, the equation
$$k_3(c|t)=\int_{-\infty}^{\infty} \; W^{'}(s-c)\; H_2(t|s) f_2(s)\;ds =0$$
has the unique solution $c\equiv \psi_{0,2}(t)$.
\\~\\ \textbf{ (i)} Let $\delta_{\psi_{0,2}}(\underline{X})=X_2-\psi_{0,2}(D)$. Then
$$R_2(\underline{\theta},\delta_{\psi_{0,2}})\leq R_2(\underline{\theta},\delta_{c_{0,2},2}), \;\;\; \forall \; \; \underline{\theta} \in \Theta_0.$$
\textbf{(ii)} Let $\psi_{1,2}:\Re \rightarrow \Re $ be such that $\psi_{1,2}(t) \geq \, (\leq) \; \psi_{0,2}(t), \; \forall \; t, \; \psi_{1,2}(t) $ is non-decreasing (non-increasing) in $t$ and $\lim_{t \to \infty}\; \psi_{1,2}(t) = c_{0,2}$.
Then 
$$R_2(\underline{\theta}, \delta_{\psi_{1,2}}) \leq R_2(\underline{\theta}, \delta_{c_{0,2},2}),\; \; \forall \; \underline{\theta} \in \Theta_0,$$
where $\delta_{\psi_{1,2}}(\underline{X})=X_2-\psi_{1,2}(D)$.\vspace{2mm}

In the following corollary we provide the Stein type improvement over the BLEE $\delta_{c_{0,2},2}(\underline{X}).$
\\~\\ \textbf{Corollary 2.2.2.} Suppose that, for any fixed $\Delta\geq 0$ and $t$, $h_2(t-\Delta\vert s)/h_2(t\vert s)$ is non-increasing (non-decreasing) in $s\in S_2$ and let $\psi_{0,2}(t)$ be as defined in Corollary 2.2.1. Further suppose that, for every $t$, the equation
$$k_4(c|t)=\int_{-\infty}^{\infty} \; W^{'}(s-c)\; h_2(t|s) f_2(s)\;ds =0$$
has the unique solution $c\equiv \psi_{2,2}(t)$. 
\\~\\ \textbf{(i)} Let $\psi_{2,2}^{*}(t)= \min\{c_{0,2},\psi_{2,2}(t)\}$ ($\psi_{2,2}^{*}(t)= \max\{c_{0,2},\psi_{2,2}(t)\}$) and $\delta_{\psi_{2,2}^{*}}(\underline{X}) = X_2- \psi_{2,2}^{*}(D)$. Then
$$R_2(\underline{\theta},\delta_{\psi_{2,2}^{*}})\leq R_2(\underline{\theta},\delta_{c_{0,2},2}), \;\;\; \forall \; \; \underline{\theta} \in \Theta_0.$$
\textbf{(ii)}  Suppose that $\psi_{3,2}:\Re \rightarrow \Re $ is such that $\psi_{3,2}(t) \geq \, (\leq) \; \psi_{2,2}(t),\,\forall \, t$ and $\psi_{3,2}(t)$ is non-decreasing (non-increasing) in $t$. Define $\psi_{3,2}^{*}(t)=\min\{c_{0,2},\psi_{3,2}(t)\}\;(\psi_{3,2}^{*}(t)$ $=\max\{c_{0,2},\psi_{3,2}(t)\})$ and $\delta_{\psi_{3,2}^{*}}(\underline{X})=X_2-\psi_{3,2}^{*}(D)$.
Then $$R_2(\underline{\theta}, \delta_{\psi_{3,2}^{*}}) \leq R_2(\underline{\theta}, \delta_{c_{0,2},2}),\; \; \forall \; \underline{\theta} \in \Theta_0.$$

It is straightforward to see that the B-Z type estimator $\delta_{\psi_{0,2}}(\cdot),$ derived in Corollary 2.2.1 (i), is the generalised Bayes estimator with respect to the non-informative prior density $\pi(\theta_1,\theta_2)=1,\;(\theta_1,\theta_2)\in\Theta_0.$ Theorems 2.1.1-2.2.1 (or Corollaries 2.1.1-2.1.2 and Corollaries 2.2.1-2.2.2) are applicable to a variety of situations studied in the literature for specific probability models, having independent marginals, and specific loss functions (e.g., Kushary and Cohen (\citeyear{MR1029476}), Misra and Singh (\citeyear{MR1366828}), Vijayasree et al. (\citeyear{MR1345425}), Misra et al. (\citeyear{MR2205815}), etc.). Theorems 2.1.1-2.2.1 (or Corollaries 2.1.1-2.1.2 and Corollaries 2.2.1-2.2.2) also extend the results of Kubokawa and Saleh (1994) to general bivariate location models.
\vspace*{3mm}

\subsection{\textbf{Applications}}
\label{sec:2.3}

\setcounter{equation}{0}
\renewcommand{\theequation}{2.3.\arabic{equation}}

\noindent
\vspace*{2mm}

In the sequel we demonstrate an application of Theorems 2.1.1-2.2.1 (or Corollaries 2.1.1-2.1.2 and Corollaries 2.2.1-2.2.2) to a situation where results of Kubokawa and Saleh (\citeyear{MR1370413}) are not applicable.
\\~\\ \textbf{Example 2.3.1.} Let $\underline{X}=(X_1,X_2)$ have the bivariate normal distribution with joint pdf given by (1.1), where, for known  $\sigma_i>0,\;i=1,2,$ and $\rho \in (-1,1),$
$$ f(z_1,z_2) =\frac{1}{2 \pi \sigma_1 \sigma_2 \sqrt{1-\rho^2}} e^{-\frac{1}{2(1-\rho^2)}\left[\frac{z_1^2}{\sigma_1^2}-2 \rho \, \frac{z_1 z_2}{\sigma_1 \sigma_2}+\frac{z_2^2}{\sigma_2^2}\right]},\; \; \; \underline{z}=(z_1,z_2)\in \Re^2.$$ 
Consider estimation of location parameter $\theta_i,\,i=1,2,$ under the squared error loss function (i.e., $W(t)=t^2,\; t\in \Re$). Here the BLEE of $\theta_i$ is $\delta_{0,i}(\underline{X})=X_i$ (i.e., $c_{0,i}=0$),$\;i=1,2$. Also, for any $s\in \Re$, and $t\in \Re$, $h_i(t\vert s)= \frac{1}{\xi_i}\,\phi \left(\frac{t-\frac{s\mu_i}{\sigma_i}}{\xi_i}\right)$ and $H_i(t\vert s)= \Phi \left(\frac{t-\frac{s\mu_i}{\sigma_i}}{\xi_i}\right),\;i=1,2$, where $\mu_1 =\rho \sigma_2-\sigma_1$, $\mu_2 =\sigma_2-\rho \sigma_1$ and $\xi_i^2=(1-\rho^2)\sigma_i^2,\;i=1,2$. For $\mu_i <(>) \;0,$ it is easy to verify that, for any fixed $\Delta\geq 0$ and $t\in \Re$, $h_i(t-\Delta\vert s)/h_i(t\vert s)$ is non-decreasing (non-increasing) in $s\in \Re,\;i=1,2$.
For any $t\in \Re$,
\begin{align*}
\psi_{0,i}(t)&= \frac{\int_{-\infty}^{\infty} s\,H_i(t|s)f_i(s)ds}{\int_{-\infty}^{\infty} H_i(t|s)f_i(s)ds}
= \frac{\int_{-\infty}^{\infty} s\,\Phi \left(\frac{t-\frac{s\mu_i}{\sigma_i}}{\xi_i}\right) \frac{1}{\sigma_i}\, \phi\left(\frac{s}{\sigma_i}\right) ds}{\int_{-\infty}^{\infty}\Phi \left(\frac{t-\frac{s\mu_i}{\sigma_i}}{\xi_i}\right) \frac{1}{\sigma_i} \,\phi\left(\frac{s}{\sigma_i}\right) ds},\;i=1,2,\\
\text{and}\qquad
\psi_{2,i}(t)&= \frac{\int_{-\infty}^{\infty} s\,h_i(t|s)f_i(s)ds}{\int_{-\infty}^{\infty} h_i(t|s)f_i(s)ds}
= \frac{\int_{-\infty}^{\infty} \frac{s}{\xi_i}\,\phi \left(\frac{t-\frac{s\mu_i}{\sigma_i}}{\xi_i}\right) \frac{1}{\sigma_i}\, \phi\left(\frac{s}{\sigma_i}\right) ds}{\int_{-\infty}^{\infty}\frac{1}{\xi_i}\phi \left(\frac{t-\frac{s\mu_i}{\sigma_i}}{\xi_i}\right) \frac{1}{\sigma_i} \,\phi\left(\frac{s}{\sigma_i}\right) ds},\;i=1,2.
\end{align*}
It is easy to verify that $\psi_{0,1}(t)=-(\beta_0-1)\,\tau\, \frac{\phi\left(\frac{t}{\tau}\right)}{\Phi\left(\frac{t}{\tau}\right)},\;t\in\Re,\; \psi_{0,2}(t)=-\beta_0\,\tau\, \frac{\phi\left(\frac{t}{\tau}\right)}{\Phi\left(\frac{t}{\tau}\right)},\;t\in\Re,\; \psi_{2,1}(t)=(\beta_0-1) \,t,\;t\in\Re,$ and $\psi_{2,2}(t)=\beta_0\, t,\;t\in\Re$, where $\tau^2=\sigma_1^2+\sigma_2^2-2\rho \sigma_1 \sigma_2$ and $\beta_0 = 1+\frac{\sigma_1 \mu_1}{\tau^2}=\frac{\sigma_2^2-\rho \sigma_1\sigma_2}{\sigma_1^2+\sigma_2^2-2\rho \sigma_1\sigma_2}=\frac{\sigma_2 \mu_2}{\tau^2}$.
\\~\\ \textbf{\underline{Estimation of $\theta_1$}:}\vspace*{2mm}

\noindent
For $\mu_1<0$, i.e., $\rho<\frac{\sigma_1}{\sigma_2}$ ($\mu_1>0$, i.e., $\rho>\frac{\sigma_1}{\sigma_2}$), we have $\beta_0<(>)1$, $\lim_{t \to \infty} \psi_{0,1}(t) = 0 = c_{0,1}$ and, $\psi_{0,1}(t)$ and $\psi_{2,1}(t)$ are non-increasing (non-decreasing) functions of $t\in\Re$. Thus, functions $\psi_{0,1}(t)$ and $\psi_{2,1}(t)$ satisfy hypotheses of Theorem 2.1.1 and Corollaries 2.1.1 and 2.1.2. For $\mu_1<0$, i.e., $\rho<\frac{\sigma_1}{\sigma_2}$ ($\mu_1>0$, i.e., $\rho>\frac{\sigma_1}{\sigma_2}$), we have
\begin{small}
$$\psi_{2,1}^{*}(t)=\max\{0,\psi_{2,1}(t)\}\;=\begin{cases}(\beta_0-1) t,& t<0 \\ 0,&t\geq 0 \end{cases}\;\;
\left(\psi_{2,1}^{*}(t)=\min\{0,\psi_{2,1}(t)\}=\begin{cases} (\beta_0-1) t,& t\leq0 \\ 0,& t> 0 \end{cases}\right).$$
\end{small}
Using Corollaries 2.1.1 (i) and 2.1.2 (i), we obtain the B-Z type and the Stein type improvements over the BLEE $\delta_{\psi_{0,1}}(\underline{X})=X_1$ as 
\begin{small}
\begin{align}
	\delta_{\psi_{0,1}}(\underline{X})&=X_1-\psi_{0,1}(D)=X_1+\frac{\sigma_1 \mu_1}{\tau} \frac{\phi\left(\frac{D}{\tau}\right)}{\Phi\left(\frac{D}{\tau}\right)}=X_1+(\beta_0-1)\tau \frac{\phi\left(\frac{D}{\tau}\right)}{\Phi\left(\frac{D}{\tau}\right)}\label{eq:2.3.1} \\
	\text{and }\; \; \delta_{\psi_{2,1}^{*}}(\underline{X})&=X_1-\psi_{2,1}^{*}(D)=\begin{cases} X_1, \qquad \qquad \qquad \; \; \; \;\; \; \; \text{if }X_1\leq X_2 \\ \beta_0 X_1+(1-\beta_0)X_2\;,  \; \; \, \text{if }X_1> X_2 \end{cases},\label{eq:2.3.2} 
\end{align}
\end{small}
respectively, where $\beta_0 = 1+\frac{\sigma_1 \mu_1}{\tau^2}=\frac{\sigma_2^2-\rho \sigma_1\sigma_2}{\sigma_1^2+\sigma_2^2-2\rho \sigma_1\sigma_2}$.\vspace*{2mm}

It is worth mentioning here that $\delta_{\psi_{0,1}}(\cdot)$ is the generalized Bayes estimator $\theta_1$ with respect to non-informative prior on $\Theta_0$ and $\delta_{\psi_{2,1}^{*}}(\cdot)$ is the restricted maximum likelihood estimator of $\theta_1$. (see Patra and Kumar \citeyear{doi:10.1080/01966324.2017.1296797}).
\\~\\  Note that, when $\mu_1=0$ (i.e., $\rho =\frac{\sigma_1}{\sigma_2}$ and $\beta_0=1$), we have $\psi_{0,1}(t)=\psi_{2,1}(t)=\psi_{2,1}^{*}(t)=c_{0,1}=0,\; \forall\; t\in \Re.$ Thus, for $\rho =\frac{\sigma_1}{\sigma_2}$, we are not able to get improvements over the BLEE using our results. Interestingly, in this case, the BLEE is also the restricted maximum likelihood estimator and the generalized Bayes estimator with respect to non-informative prior on $\Theta_0$.
\\~\\ From the above discussion we conclude that, for $\rho \neq \frac{\sigma_1}{\sigma_2}$, the generalized Bayes estimator $\delta_{\psi_{0,1}}(\cdot)$ and the restricted MLE $\delta_{\psi_{2,1}^{*}}(\cdot)$ dominate the BLEE $\delta_{0,1}(\underline{X})=X_1$.
\\~\\ Now, we will illustrate an application of Corollary 2.1.1 (ii). Define 
$$\psi_{1,1,\alpha}(t)=\tau \,(1-\alpha)\, \frac{\phi\left(\frac{t}{\tau}\right)}{\Phi\left(\frac{t}{\tau}\right)},\; t\in \Re,\;\alpha\in \Re.$$
For $\mu_1<0$ and $\beta_0\leq \alpha<1$ ($\mu_1>0$ and $1<\alpha\leq \beta_0$), note that $\beta_0<1$ ($\beta_0>1$), $\psi_{1,1,\alpha}(t)$ is a non-increasing (non-decreasing) functions of $t\in \Re$, $\lim_{t\to \infty} \psi_{1,1,\alpha}(t)=0=c_{0,1}$ and $\psi_{1,1,\alpha}(t)\leq\; (\geq)\; \psi_{0,1}(t),\,\forall\, t\in \Re.$ Let
$$\delta_{\psi_{1,1,\alpha}}(\underline{X})=X_1-\psi_{1,1,\alpha}(D)= X_1-\tau \,(1-\alpha)\, \frac{\phi\left(\frac{D}{\tau}\right)}{\Phi\left(\frac{D}{\tau}\right)},\; \alpha\in \Re.$$
Using Corollary 2.1.1 (ii) it follows that, for $\mu_1\!<\! (>)\,0$ (i.e., $\rho <\!(>)\,\frac{\sigma_1}{\sigma_2}$), the estimators $\{\delta_{\psi_{1,1,\alpha}}\,:\, \beta_0\leq \alpha<1\}$ ($\{\delta_{\psi_{1,1,\alpha}}\,:\, 1<\alpha\leq \beta_0\}$) dominate the BLEE $\delta_{0,1}(\underline{X})=X_1.$
\\To see an application of Corollary 2.1.2 (ii), let 
$$\psi_{3,1,\alpha}(t)=\begin{cases} (\alpha-1)\,t, \;\; \; t<0\\ (\beta_0-1)\,t,\; \; t\geq 0 \end{cases} , \; \; \; \alpha\in \Re.$$ 
For $\mu_1<0$ and $\beta_0\leq \alpha < 1$ ($\mu_1>0$ and $1<\alpha\leq \beta_0$), note that $\beta_0<1$ ($\beta_0>1$), $\psi_{3,1,\alpha}(t)\leq \,(\geq)\;\psi_{2,1}(t)=(\beta_0-1)\,t , \; \forall \; t\in \Re$, and $\psi_{3,1,\alpha}(t)$ is a non-increasing (non-decreasing) function of $t\in \Re$. Let
\begin{align*}
\psi_{3,1,\alpha}^{*}(t)&=\max\{0,\psi_{3,1,\alpha}(t)\}\;(\min\{0,\psi_{3,1,\alpha}(t)\})
=\begin{cases} (\alpha-1)\,t,\; \; t<0\\ 0,\;\quad \; \qquad\; t\geq 0 \end{cases}\\
\text{and}\quad \delta_{\psi_{3,1,\alpha}^{*}}(\underline{X})&=X_1-\psi_{3,1,\alpha}^{*}(D)=\begin{cases} \alpha X_1+(\alpha-1)X_2\,,\; \; X_2<X_1\\ X_1 ,\qquad \qquad \qquad \quad X_2\geq X_1 \end{cases},\; \; \alpha\in\Re.
\end{align*} 
Using Corollary 2.1.2 (ii), it follows that, for $\rho<(>)\, \frac{\sigma_1}{\sigma_2}$, the estimators $\{\delta_{\psi_{3,1,\alpha}^{*}}\,:\, \beta_0\leq \alpha <1\}$ ($\{\delta_{\psi_{3,1,\alpha}^{*}}\,:\, 1<\alpha\leq \beta_0\}$) dominate the BLEE $\delta_{0,1}(\underline{X})=X_1.$
\\~\\ \textbf{\underline{Estimation of $\theta_2$}:}\vspace*{2mm}

\noindent
For $\mu_2<0,$ i.e., $\rho >\frac{\sigma_2}{\sigma_1}$ ($\mu_2>0,$ i.e., $\rho <\frac{\sigma_2}{\sigma_1}$), we have $\beta_0<(>)0$, $\lim_{t \to \infty} \psi_{0,2}(t) = 0 = c_{0,2}$ and, $\psi_{0,2}(t)$ and $\psi_{2,2}(t)$ are non-increasing (non-decreasing) functions of $t\in\Re$. Let
\begin{small}
$$\psi_{2,2}^{*}(t)=\max\{0,\psi_{2,2}(t)\}=\begin{cases} \beta_0 t,&\text{ if  } t\leq0 \\ 0,&\text{ if  } t> 0 \end{cases}\;\left(\psi_{2,2}^{*}(t)=\min\{0,\psi_{2,2}(t)\}=\begin{cases} \beta_0 t,&\text{ if  } t\leq0 \\ 0,&\text{ if  } t> 0 \end{cases}\right).$$
\end{small}
Applications of Corollaries 2.2.1 (i) and 2.2.2 (i), yield the B-Z type and the Stein type improvements over the BLEE $\delta_{0,2}(\underline{X})=X_2$ as 
\begin{align*}
\delta_{\psi_{0,2}}(\underline{X})&=X_2-\psi_{0,2}(D)=X_2+\frac{\sigma_2 \mu_2}{\tau} \frac{\phi\left(\frac{D}{\tau}\right)}{\Phi\left(\frac{D}{\tau}\right)}=X_2+\beta_0\tau \frac{\phi\left(\frac{D}{\tau}\right)}{\Phi\left(\frac{D}{\tau}\right)}\\
\text{and }\; \; \delta_{\psi_{2,2}^{*}}(\underline{X})&=X_2-\psi_{2,2}^{*}(D)=\begin{cases}\beta_0 X_1+(1-\beta_0)X_2, \;\;\;\; \text{if }X_2\leq X_1 \\ X_2\;,  \qquad \qquad \qquad \; \; \; \;\; \; \; \text{if }X_2> X_1 \end{cases},
\end{align*}
respectively. Note that $\delta_{\psi_{0,2}}(\cdot)$ is the generalized Bayes estimators of $\theta_2$ under the non-informative prior on $\Theta_0$ and $\delta_{\psi_{2,2}^{*}}(\cdot)$ is the restricted MLE of $\theta_2$.
\\~\\ For $\mu_2=0$ (i.e., $\rho =\frac{\sigma_2}{\sigma_1}$), we have $\beta_0=0$ and $\psi_{0,2}(t)=\psi_{2,2}(t)=\psi_{2,2}^{*}(t)=0,\; \forall\; t\in \Re.$ Thus, for $\mu_2=0$, our results do not provide improvements over the BLEE $\delta_{\psi_{0,2}}(\underline{X})=X_2$. In this case, the BLEE is also the restricted maximum likelihood estimator and the generalized Bayes estimator with respect to non-informative prior on $\Theta_0$.
\\~\\ From the above discussion we conclude that, for $\rho \neq \frac{\sigma_2}{\sigma_1}$, the generalized Bayes estimator $\delta_{\psi_{0,2}}(\cdot)$ and the restricted MLE $\delta_{\psi_{2,2}^{*}}(\cdot)$ dominate the BLEE $\delta_{0,2}(\underline{X})=X_2$.
\\To see an application of Corollary 2.2.1 (ii), define 
$$\psi_{1,2,\alpha}(t)=-\alpha\,\tau \, \frac{\phi\left(\frac{t}{\tau}\right)}{\Phi\left(\frac{t}{\tau}\right)},\; t\in \Re,\;\alpha\in \Re.$$
For $\mu_2<0$ and $\beta_0\leq \alpha<0$ ($\mu_2>0$ and $0<\alpha\leq \beta_0$), note that $\beta_0<(>)0$, $\psi_{1,2,\alpha}(t)$ is a non-increasing (non-decreasing) function of $t\in \Re$, $\lim_{t\to \infty} \psi_{1,2,\alpha}(t)=0=c_{0,2}$ and $\psi_{1,2,\alpha}(t)\leq\; (\geq)\; \psi_{0,2}(t)=\beta_0 \tau \,\frac{\phi\left(\frac{t}{\tau}\right)}{\Phi\left(\frac{t}{\tau}\right)},\; \forall\; t\in \Re.$ Define 
$$\delta_{\psi_{1,2,\alpha}}(\underline{X})=X_2-\psi_{2,2,\alpha}(D)= X_2+\alpha\,\tau \, \frac{\phi\left(\frac{D}{\tau}\right)}{\Phi\left(\frac{D}{\tau}\right)},\; \alpha\in \Re.$$
Using Corollary 2.2.1 (ii), for $\rho>(<)\,\frac{\sigma_2}{\sigma_1}$, it follows that the estimators $\{\delta_{\psi_{1,2,\alpha}}\,:\, \beta_0\leq \alpha<0\}$ ($\{\delta_{\psi_{1,2,\alpha}}\,:\, 0<\alpha\leq \beta_0\}$) dominate the BLEE $\delta_{0,2}(\underline{X})=X_2.$
\\Now consider an application of Corollary 3.2.2 (ii). Define 
$$\psi_{3,2,\alpha}(t)=\begin{cases} \alpha\,t,\;\; \; t<0\\ \beta_0\,t,\; \; t\geq 0 \end{cases} , \; \; \; \alpha\in \Re.$$ 
For $\mu_2<0$ and $\beta_0\leq \alpha<0$ ($\mu_2>0$ and $0<\alpha\leq \beta_0$) note that, $\beta_0<(>)\,0$, $\psi_{3,2,\alpha}(t)$ is non-increasing (non-decreasing) in $t\in \Re$ and  $\psi_{3,2,\alpha}(t)\leq \,(\geq)\;\psi_{2,2}(t)=\beta_0 \,t , \; \forall \; t\in \Re$. Let $$\psi_{3,2,\alpha}^{*}(t)=\max\{0,\psi_{3,2,\alpha}(t)\}\;(\min\{0,\psi_{3,2,\alpha}(t)\})=\begin{cases} \alpha\,t, \quad t<0\\ 0, \; \quad\; t\geq 0 \end{cases}$$ 
\\and $\qquad\,\,\quad\delta_{\psi_{3,2,\alpha}^{*}}(\underline{X})=X_2-\psi_{3,2,\alpha}^{*}(D)=\begin{cases} \alpha X_1+(\alpha-1)X_2,\; \; X_2<X_1\\ X_2, \;\qquad \; \qquad \; \qquad X_2\geq X_1 \end{cases}.$ 
\\Using Corollary 2.2.2 (ii), it follows that, for $\rho>(<)\, \frac{\sigma_2}{\sigma_1}$, the class of the estimators $\{\delta_{\psi_{3,2,\alpha}^{*}}\,:\, \beta_0\leq \alpha< 0\}$ ($\{\delta_{\psi_{3,2,\alpha}^{*}}\,:\, 0< \alpha \leq \beta_0\}$) dominate the BLEE $\delta_{0,2}(\underline{X})=X_2.$

\subsection{\textbf{Simulation Study For Estimation of Location Parameter $\theta_1$}}
\label{sec:2.4}

\noindent
\vspace*{2mm}

In Example 2.2.1, under the squared error loss function, we have considered estimation of the smaller mean $\theta_1$ of a bivariate normal distribution with unknown order restricted

\FloatBarrier
\begin{figure}[h!]
\begin{subfigure}{0.48\textwidth}
	%\centering
	% include first image
	\includegraphics[width=68mm,scale=1.2]{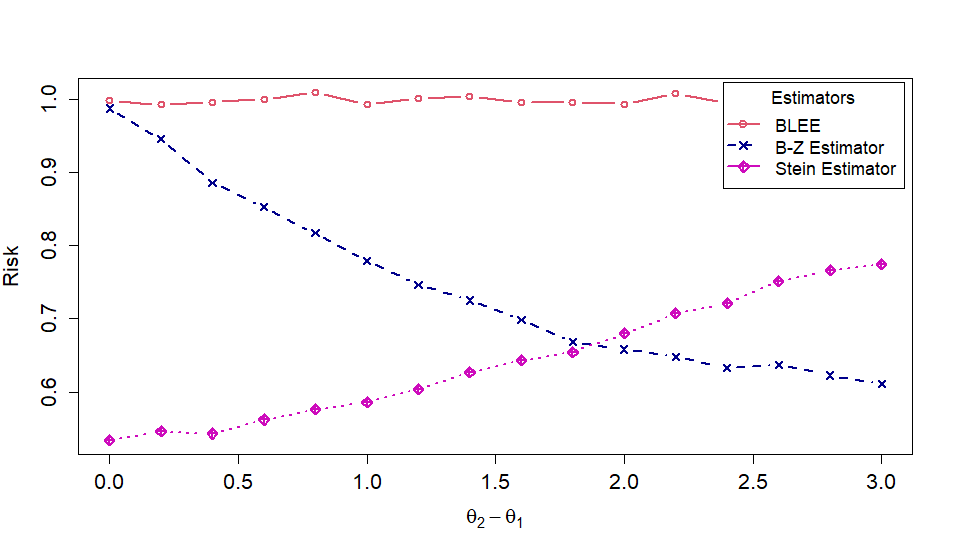} 
	\caption{$\sigma_1=0.2$, $\sigma_2=0.2$ and $\rho=-0.9$.} 
	\label{fig7:a} 
\end{subfigure}
\begin{subfigure}{0.48\textwidth}
	%	\centering
	% include third image
	\includegraphics[width=68mm,scale=1.2]{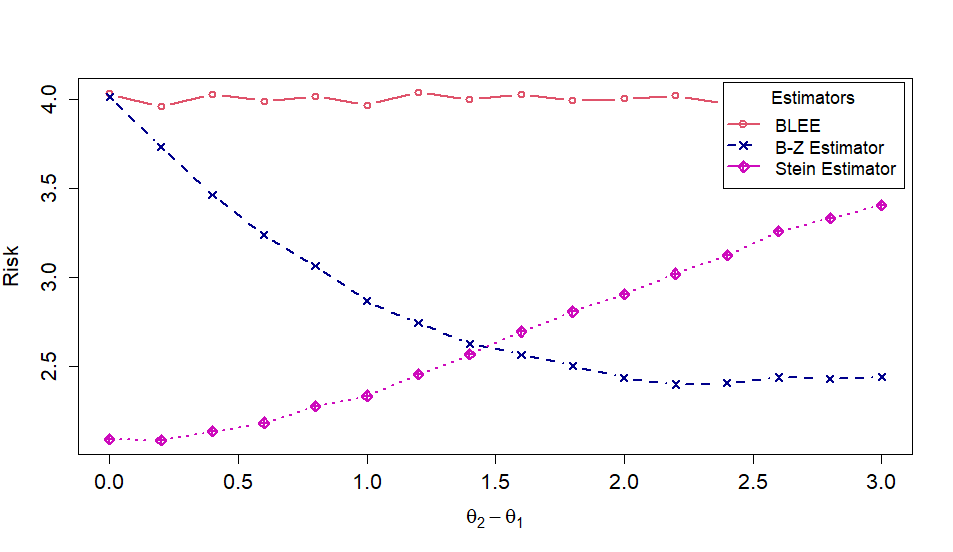} 
	\caption{$\sigma_1=2$, $\sigma_2=0.5$ and $\rho=-0.5$.} 
	\label{fig7:b} 
\end{subfigure}
\\	\begin{subfigure}{0.48\textwidth}
	%	\centering
	% include fourth image
	
	\includegraphics[width=68mm,scale=1.2]{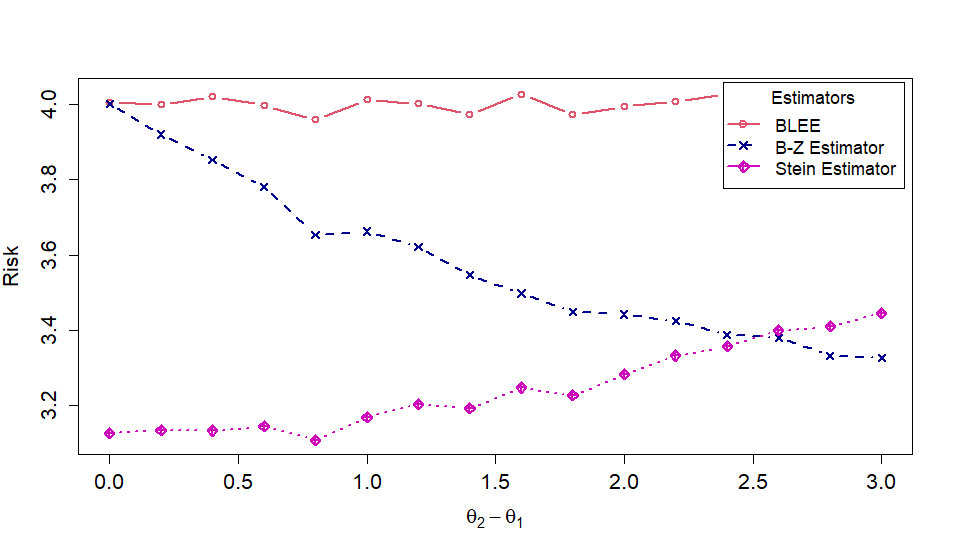} 
	
	\caption{$\sigma_1=2$, $\sigma_2=3$ and $\rho=-0.2$.} 
	\label{fig7:c} 
\end{subfigure}
\begin{subfigure}{0.48\textwidth}
	%	\centering
	% include third image
	\includegraphics[width=68mm,scale=1.2]{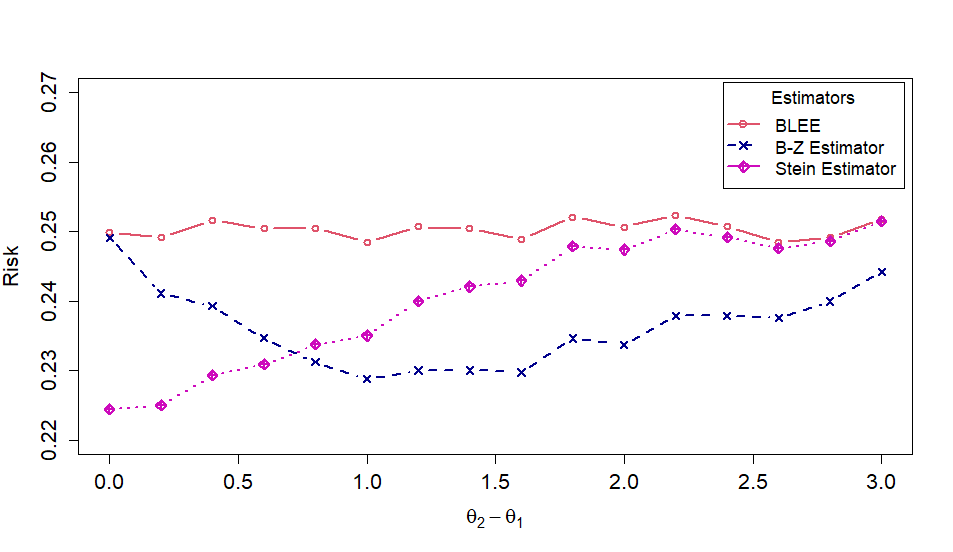} 
	\caption{ $\sigma_1=0.5$, $\sigma_2=1$ and $\rho=0$.} 
	\label{fig7:d}  
\end{subfigure}
\\	\begin{subfigure}{0.48\textwidth}
	\centering
	% include fourth image
	
	\includegraphics[width=68mm,scale=1.2]{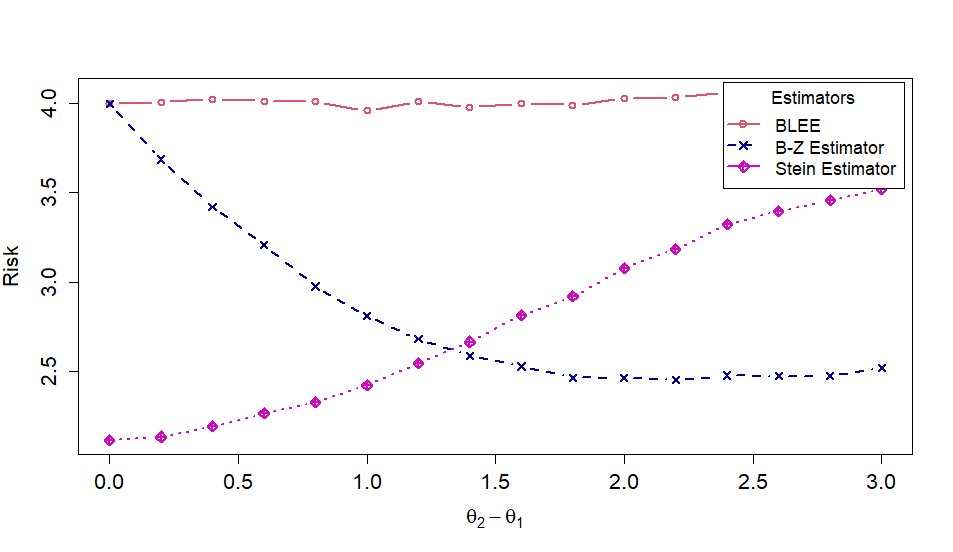} 
	\caption{$\sigma_1=2$, $\sigma_2=0.5$ and $\rho=0$.} 
	\label{fig7:e} 
\end{subfigure}
\begin{subfigure}{0.48\textwidth}
	\centering
	% include third image
	\includegraphics[width=68mm,scale=1.2]{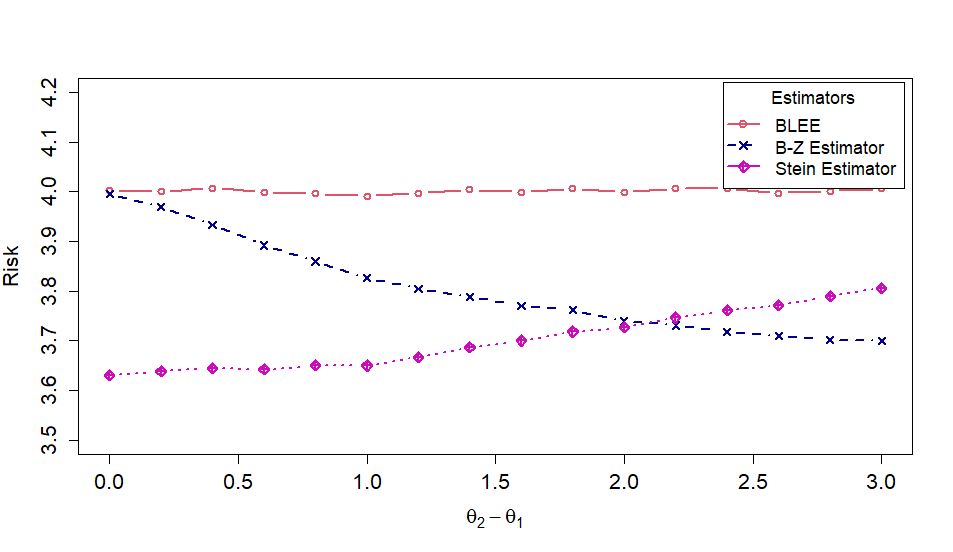} 
	\caption{ $\sigma_1=2$, $\sigma_2=3$ and $\rho=0.2$.} 
	\label{fig7:f}  
\end{subfigure}
\\	\begin{subfigure}{.48\textwidth}
	\centering
	% include third image
	\includegraphics[width=68mm,scale=1.2]{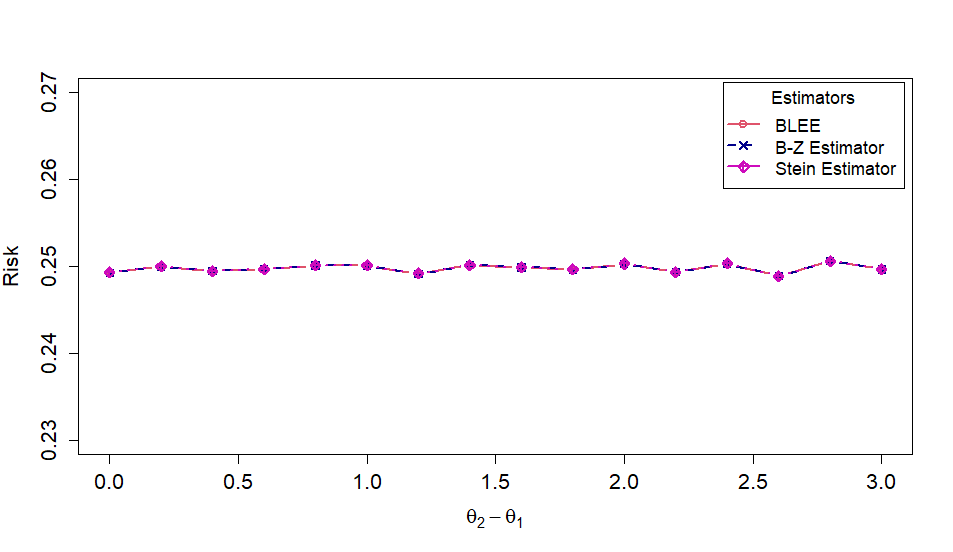} 
	\caption{ $\sigma_1=0.5$, $\sigma_2=1$ and $\rho=0.5$.} 
	\label{fig7:g}  
\end{subfigure}
\begin{subfigure}{.48\textwidth}
	\centering
	% include fourth image
	
	\includegraphics[width=68mm,scale=1.2]{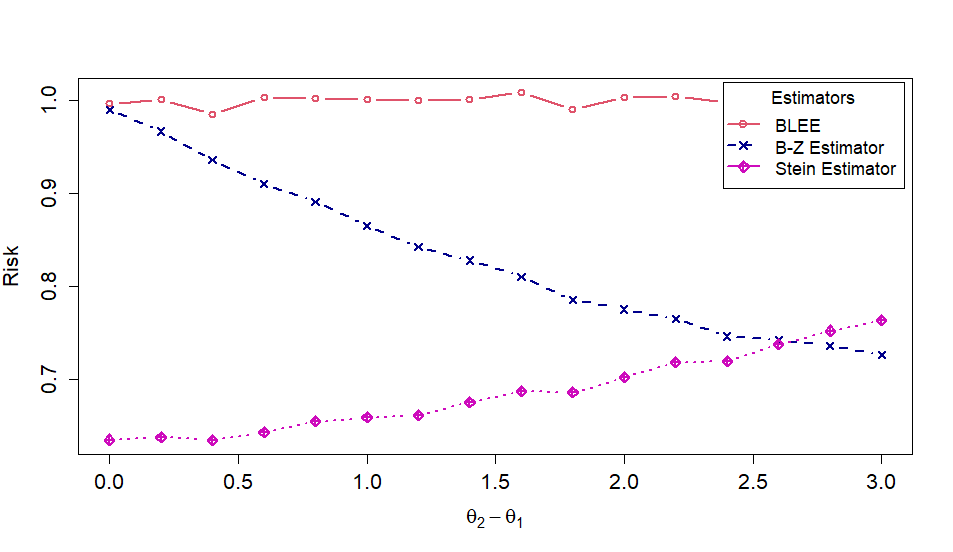} 
	\caption{$\sigma_1=1$, $\sigma_2=5$ and $\rho=0.9$.} 
	\label{fig7:h} 
\end{subfigure}

\caption{Risk plots of $\delta_{0,1}$ (BLEE), $\delta_{\psi_{0,1}}$ (B-Z type estimator) and $\delta_{\psi^{*}_{2,1}}$ (Stein type estimator) estimators against the values of $\theta_2-\theta_1$}
\label{fig7}
\end{figure}
\FloatBarrier

\newpage

\noindent
 means (i.e., $\theta_1\leq \theta_2$), known variances ($\sigma_1^2$ and $\sigma_2^2$) and known correlation coefficient ($\rho$), and obtained improvements over the BLEE $\delta_{0,1}(\underline{X})=X_1$. To further evaluate the performances of various estimators under the squared error loss function, in this section, we compare the risk performances of estimators BLEE $\delta_{0,1}(\underline{X})=X_1$, the B-Z estimator $\delta_{\psi_{0,1}}$ and the Stein (1964) type estimator $\delta_{\psi^{*}_{2,1}}$ (as defined in \eqref{eq:2.3.1} and \eqref{eq:2.3.2}), numerically, through Monte Carlo simulations. The simulated risks of the BLEE, the B-Z estimator and the Stein estimator (restricted MLE) have been computed based on 10000 simulations from relevant distributions. Note that the B-Z estimator is the generalized Bayes estimator of $\theta_1$ and the Stein estimator is the restricted MLE of $\theta_1$ under $\underline{\theta}\in \Theta_0$.

The simulated values of risks of various estimators are plotted in Figure \ref{fig7}. The following observations are evident from Figure \ref{fig7}:

\noindent
(i) The risk function values of the B-Z type and the Stein type estimators are nowhere larger than the risk function values of the BLEE, which is in conformity with theoretical findings of Example 2.2.1.
\\(ii) There is no clear cut winner between the B-Z type estimator $\delta_{\psi_{0,1}}$ and the Stain type estimator $\delta_{\psi^{*}_{2,1}}$. The Stein type estimator performs better than the B-Z type estimator, for small values of $\theta_2-\theta_1$, and the B-Z type estimator dominates the Stein type estimator for the large values of $\theta_2-\theta_1$.

\section{\textbf{Improving the Best Scale Equivariant Estimators (BSEEs)}}
\label{sec:3}

\setcounter{equation}{0}
\renewcommand{\theequation}{3.\arabic{equation}}

In this section, we consider the bivariate scale model \eqref{eq:1.2}, and deal with the problem of estimating scale parameters $\theta_i,\;i=1,2,$ when it is known apriori that $\underline{\theta}\in\Theta_0=\{(x,y)\in \Re_{++}^2\,:\,x\leq y\}$. The following notations will be used throughout this section.
Let $Z_i=\frac{X_i}{\theta_i}$, $i=1,2,$ \,$\underline{Z}=(Z_1,Z_2)$ and  $Z=\frac{Z_2}{Z_1}$. The pdf of $\underline{Z}=(Z_1,Z_2)$ is $f(z_1,z_2)$, $(z_1,z_2)\in \Re^2$. Let $S_i$ denote the support of random variable $Z_i,\;i=1,2.$ Under the above notations, assume that $\{(z_1,z_2)\in\Re^2:\;f(z_1,z_2)>0\}\subseteq \Re_{++}^2$, so that $S_1\subseteq\Re_{++}$ and $S_2\subseteq \Re_{++}$. Let $f_i$ denote the pdf of $Z_i,\;i=1,2$, so that $f_1(s)=\int\limits_{0}^{\infty} f(s,t) \,dt,\; s\in \Re_{++} \;\text{ and }\; f_2(s)=\int\limits_{0}^{\infty} f(t,s) \,dt,\;\; s\in \Re_{++}.$

For any $s\in S_i$, let $Z_s^{(i)}$ denote a random variable having the same distribution as conditional distribution of $Z$ given $Z_i=s,\;i=1,2$. Then, the pdf and the df of $Z_s^{(1)}\,(s\in S_1)$ are given by
\begin{equation*}
h_1(t|s)=s\frac{f(s,st)}{f_1(s)}, \;\; t\in\Re_{++}, \text{ and }\; H_1(t|s)=\int_{0}^{t}h_1(z|s)\,dz, \;\; t\in\Re_{++},
\end{equation*}
respectively and the pdf and the df of $Z_s^{(2)}\,(s\in S_2)$ are given by
\begin{equation*}
h_2(t|s)=\frac{s}{t^2}\frac{f\left(\frac{s}{t},s\right)}{f_2(s)}, \;\; t\in\Re_{++}, \text{ and }\; H_2(t|s)=\int_{0}^{t}	h_2(z|s)\,dz, \;\; t\in\Re_{++},
\end{equation*}
respectively.
\\~\\For the scale model (1.2), consider estimation of scale parameter $\theta_i$ under the loss function 
\begin{equation}\label{eq:3.1}
L_i(\underline{\theta},a)=W\left(\frac{a}{\theta_i}\right),  \;\underline{\theta}\in\Theta,\; a\in\mathcal{A} =\Re_{++},\;i=1,2,
\end{equation}    
where $W:\Re\rightarrow [0,\infty)$ is a specified non-negative function. Throughout, we make the following assumptions on the function $W(\cdot)$:
\\~\\ \textbf{\boldmath$A_3$:} $W:\Re \rightarrow [0,\infty)$ is such that $W(1)=0$, $W(t)$ is decreasing on $(-\infty,1)$ and increasing on $(1,\infty)$. Further $W^{'}(t)$ is non-decreasing on the set $D_0$ (the set of points at which $W(\cdot)$ is differentiable).
\\~\\ \textbf{\boldmath$A_4$:} The equation $E[Z_i\,W^{'}(cZ_i)]=0$ has the unique solution, say $c=c_{0,i},\;i=1,2.$ \vspace{2mm}

Under the unrestricted case ($\Theta=\Re_{++}$), the problem of estimating $\theta_i$, under the loss function (3.2) is invariant under the multiplicative group of transformations $\mathcal{G}_0\!=\!\{g_{b_1,b_2}\!:\!\,(b_1,b_2)\!\in\Re_{++}^2\},$ where $g_{b_1,b_2}(x_1,x_2)=(b_1x_1,b_2x_2),\; (x_1,x_2)\in\Re^2,\;(b_1,b_2)\in\Re_{++}^2$, and the best scale equivariant estimator of $\theta_i$ is
$\delta_{c_{0,i},i}(\underline{X})=c_{0,i}X_i,\;i=1,2,$
where $c_{0,i}$ is the unique solution of the equation 
$ \int_{0}^{\infty}\,s \,W'(cs) \,f_i(s)ds=0,\;i=1,2$.

Under the restricted parameter space $\Theta_0$,
the problem of estimating $\theta_i$, under the loss function \eqref{eq:3.1}, is invariant under the group of transformations $\mathcal{G}\!=\!\{g_b:\,b\in(0,\infty)\},$ where $g_b(x_1,x_2)\!=\!(b\,x_1,b\,x_2),\; (x_1,x_2)\in\Re^2,\;b\in(0,\infty)$. Any scale equivariant estimator of $\theta_i$ has the form
\begin{equation}\label{eq:3.2}
\delta_{\psi_i}(\underline{X})=\psi_i(D)X_i,
\end{equation}
for some function $\psi_i:\,\Re_{++}\rightarrow \Re\,,\;i=1,2,$ where $D=\frac{X_2}{X_1}$. The risk function
\begin{equation}\label{eq:3.3}
R_i(\underline{\theta},\delta_{\psi_i})=E_{\underline{\theta}}[L_i(\underline{\theta},\delta_{\psi_i}(\underline{X}))],\; \, \underline{\theta}\in\Theta_0,
\end{equation}
of any scale equivariant estimator $\delta_{\psi_i}$ of $\theta_i,\;i=1,2,$ depends on $\underline{\theta}\in\Theta_0$ only through $\lambda=\frac{\theta_2}{\theta_1}\in[1,\infty)$. \vspace*{1.5mm}

\noindent
The following dual of Lemma 2.2 will be useful in proving the results of this section.
\\~\\ \textbf{Lemma 3.1.} If, for any fixed $\Delta\geq 1$ and $t$, $h_i\left(\frac{t}{\Delta}\vert s\right)/h_i(t\vert s)$ is non-decreasing (non-increasing) in $s\in S_i$,
then $H_i(\frac{t}{\Delta}\vert s)/H_i(t\vert s)$ is non-decreasing (non-increasing) in $s\in S_i$ and $h_i(t \vert s)/H_i(t\vert s)$ is also non-increasing (non-decreasing) in $s\in S_i,\;i=1,2.$
\vspace*{2mm}

In Subsection 3.1 (3.2), we consider the equivariant estimation of scale parameter $\theta_1$ ($\theta_2$) under the loss function $L_1\,(L_2)$ defined by \eqref{eq:3.1}, when it is known apriori that $\underline{\theta}\in\Theta_0$. In Subsection 3.3, we provide an application of our results to a bivariate gamma distribution, not studied before in the literature. In Subsection 3.4, we report a simulation study on comparison of various competing estimators for smaller scale parameter in the Cheriyan and Ramabhadran's bivariate gamma distribution.

\subsection{\textbf{Improvements Over the BSEE of $\theta_1$}}
\label{sec:3.1}
\setcounter{equation}{0}
\renewcommand{\theequation}{3.1.\arabic{equation}}

\noindent
\vspace*{2mm}

The following theorem provides a class of estimators that improve upon the BSEE $\delta_{c_{0,1},1}(\underline{X})=c_{0,1} X_1$, where $c_{0,1}$ is the unique solution of the equation $\int_{-\infty}^{\infty} z\,W^{'}(cz)f_1(z) \,dz=0.$  
\\~\\\textbf{Theorem 3.1.1.} Suppose that, for any fixed $\Delta\geq 1$ and $t$, $H_1(\frac{t}{\Delta}\vert s)/H_1(t\vert s)$ is non-decreasing (non-increasing) in $s\in S_1$. Consider a scale equivariant estimator $\delta_{\psi_1}(\underline{X})=\psi_1(D)X_1$ for estimating $\theta_1$, where $\lim_{t\to\infty} \psi_1(t)=c_{0,1}$, $\psi_1(t)$ is a non-decreasing (non-increasing) function of $t$ and $\int_{0}^{\infty}s\, W^{'}(\psi_1(t)s)\; H_1(t\vert s)\,f_1(s)ds \, \geq \,(\leq) \,0,\; \forall\;t$. Then, $\forall\; \underline{\theta}\in \Theta_0$, the estimator $\delta_{\psi_1}(\underline{X})$ dominates the BSEE $\delta_{c_{0,1},1}(\underline{X})=c_{0,1}X_1$.
\begin{proof}
For $\underline{\theta}\in\Theta_0$ and $\lambda=\frac{\theta_2}{\theta_1}$, so that $\lambda\geq 1$, the risk difference can be written as
\begin{align*}\Delta_1(\lambda)
	& = E\left[W\left(\frac{c_{0,1}\,X_1}{\theta_1}\right)\right]- E\left[W\left(\frac{\psi_1(D)X_1}{\theta_1}\right)\right]\\
	&= \int_{0}^{\infty} \psi_1^{'}(t)\; E\left[Z_1\,W^{'}(\psi_1(t)Z_1)\,I_{\left(0,\frac{t}{\lambda} \right]}(Z)\right]\,dt.
\end{align*}

In light of the hypotheses of the theorem, it is enough to prove that, for every fixed $t$,
\begin{small}
	\begin{equation}  
		E_{\underline{\theta}}\left[Z_1\,W^{'}(\psi_1(t)Z_1)\; I_{\left(0,\frac{t}{\lambda} \right]}(Z)\;\right]=\int_{0}^{\infty} s\,W^{'}(\psi_1(t)s)\;H_1\!\left(\frac{t}{\lambda}|s\right)f_1(s)ds\; \geq ( \leq)\; 0.
	\end{equation}
\end{small}
To prove the inequality (3.1.1), let us fix $t$ and $\lambda\geq 1$. Let $s_0 =\frac{1}{\psi_1(t)},\; M_1(s)=H_1(t|s),\; M_2(s)$ $=H_1\left(\frac{t}{\lambda}|s\right)$ and $M(s)=s\, W^{'}\left(\frac{s}{s_0}\right) \; f_1(s),\; s\in S_1$. Since $W^{'}(t)\leq 0,\; \forall\; t<1$ and $W^{'}(t)$ $\geq 0,\; \forall\; t>1$, we have $M(s)\leq 0,\; \forall \; 0<s<s_0$ and $M(s)\geq 0,\; \forall \; s> s_0>0.$
Now using hypotheses of the theorem and the Lemma 2.1, we conclude that
\begin{small}
	\begin{align*}
		0\leq(\geq)\;H_1\left(\frac{t}{\lambda}|\frac{1}{\psi_1(t)}\right) \int_{-\infty}^{\infty}\;s\, W^{'}\!(\psi_1(t)s)\; H_1(t|s) f_1(s) ds \qquad \qquad \qquad \qquad \qquad \qquad \qquad \quad \\ \qquad \qquad \qquad \qquad \qquad \qquad \leq\; (\geq)\; H_1\left(t|\frac{1}{\psi_1(t)}\right) \int_{-\infty}^{\infty}\;s\, W^{'}\!(\psi_1(t)s)\; H_1\!\left(\frac{t}{\lambda}|s\right) f_1(s) ds .
	\end{align*}
\end{small}
which, in turn, implies (3.1.1).
\end{proof}

\noindent
The following corollary gives us the B-Z type improvements over the BSEE $\delta_{c_{0,1},1}(\underline{X})$.
\\~\\	\textbf{Corollary 3.1.1.} Suppose that assumptions of Theorem 3.1.1 hold. Further suppose that, for every fixed $t$,  the equation
$$l_1(c\vert t)=\int_{0}^{\infty} \;s\, W^{'}(cs)\; H_1(t\vert s)\,f_1(s)ds =0$$
has the unique solution $c\equiv \psi_{0,1}(t)$.
\\~\\ \textbf{(i)} Then the estimator $\delta_{\psi_{0,1}}(\underline{X})=\psi_{0,1}(D)X_1$ dominates the BSEE $\delta_{c_{0,1},1}(\underline{X})$, $\forall\;\underline{\theta} \in \Theta_0$. 
\\~\\\textbf{(ii)} Suppose that $\psi_{1,1}:\Re_{++} \rightarrow \Re $ is such that $\psi_{1,1}(t) \geq \; (\leq) \; \psi_{0,1}(t), \; \; \forall \;t, \; \psi_{1,1}(t) $ is non-decreasing (non-increasing) in $t$ and $\lim_{t \to \infty}\; \psi_{1,1}(t) = c_{0,1}$.
Then the estimator $\delta_{\psi_{1,1}}(\underline{X})=\psi_{1,1}(D)X_1$ dominates the BSEE $\delta_{c_{0,1},1}(\underline{X})$, $\forall\;\underline{\theta} \in \Theta_0.$
\begin{proof}
It is sufficient to prove that $\psi_{0,1}(t)$ satisfies conditions of Theorem 3.1.1. To prove that $\psi_{0,1}(t)$ is an non-decreasing (non-increasing) function of $t$, suppose that, there exist positive numbers $t_1$ and $t_2$ such that $t_1<t_2$ and $\psi_{0,1}(t_1)\neq \psi_{0,1}(t_2).$ Then $l_1(\psi_{0,1}(t_1)\vert t_1)=0$.
Since $W^{'}(t)$ is an non-decreasing function of $t\in\Re$, it follows that $l_1(c|t_2)$ is a non-decreasing function of $ c $ and $\psi_{0,1}(t_2)$ is the unique solution of $l_1(c|t_2)=0$. Let $ s_0 = \frac{1}{\psi_{0,1}(t_1)}, \; M(s)=s\,W^{'}\left(\frac{s}{s_0}\right)\,f_1(s),\; M_1(s)=H_1(t_2|s)$ and $M_2(s)=H_1(t_1|s),\; s\in S_1$. Then, using hypotheses of the corollary and the Lemma 2.1, we get 
$$l_1(\psi_{0,1}(t_1)|t_2)=\int_{0}^{\infty}\;s \,W^{'}(\psi_{0,1}(t_1)s)\, H_1(t_2|s) f_1(s) ds<(>)\,0,$$ 
as $l_1(c|t_2)=0$ has the unique solution $c=\psi_{0,1}(t_2)$ and $\psi_{0,1}(t_1)\neq\psi_{0,1}(t_2)$. Since $l_1(c|t_2)$ is a non-decreasing function of c and $l_1(\psi_{0,1}(t_2)|t_2)=0$, this implies that $\psi_{0,1}(t_1)<(>)$ $\psi_{0,1}(t_2)$.
\\~\\ Also, $l_1(\psi_{0,1}(t)\vert t)=0$ and the assumption $A_2$ ensures that $\lim_{t\to\infty} \psi_{0,1}(t)=c_{0,1}$. Hence the assertion follows.
\\~\\ The proof of part \textbf{(ii)} is immediate from Theorem 3.1.1, since $l_1(c|t)$ is a non-decreasing function of $c\in\Re$, $\forall\;t$.
\end{proof}

\noindent
The following corollary gives us the Stein type improvements over the BSEE $\delta_{c_{0,1},1}(\underline{X}).$
\\~\\ \textbf{Corollary 3.1.2.} Suppose that, for any fixed $\Delta\geq 1$ and $t$, $h_1\!\left(\frac{t}{\Delta}\vert s\right)/h_1(t\vert s)$ is non-decreasing (non-increasing) in $s\in S_1$. In addition suppose that, for any $t$, the equation
$$l_2(c\vert t)=\int_{0}^{\infty} \;s\, W^{'}(c\,s)\; h_1(t\vert s)\,f_1(s)\;ds =0$$
has the unique solution $c\equiv \psi_{2,1}(t)$.
\\~\\ \textbf{(i)} Let $\psi_{2,1}^{*}(t)= \min\{c_{0,1},\psi_{2,1}(t)\}$ ($\psi_{2,1}^{*}(t)= \max\{c_{0,1},\psi_{2,1}(t)\}$) and $\delta_{\psi_{2,1}^{*}}(\underline{X}) =  \psi_{2,1}^{*}(D)X_1$. Then, $\forall \; \underline{\theta} \in \Theta_0$, the estimator $\delta_{\psi_{2,1}^{*}}(\underline{X})$ dominates the BSEE $\delta_{c_{0,1},1}(\underline{X})=c_{0,1}X_1$.
\\~\\\textbf{(ii)} Let $\psi_{3,1}:\Re_{++} \rightarrow \Re $ be such that $\psi_{3,1}(t) \geq (\leq) \, \psi_{2,1}(t), \, \forall \,t$ and $\psi_{3,1}(t) $ is non-decreasing (non-increasing) in $t$. 
Define $\psi_{3,1}^{*}(t)=\min\{c_{0,1},\psi_{3,1}(t)\}\;(\psi_{3,1}^{*}(t)=\max\{c_{0,1},\psi_{3,1}(t)\})$. Then, $\forall \; \underline{\theta} \in \Theta_0$, the estimator $\delta_{\psi_{3,1}^{*}}(\underline{X})=\psi_{3,1}^{*}(D)X_1$ dominates the BSEE $\delta_{c_{0,1},1}(\underline{X})$.
\begin{proof}
It suffices to show that $\psi_{2,1}^{*}(\cdot)$ satisfies conditions of Theorem 3.1.1. On using arguments similar to the ones used in the proof of Corollary 3.1.1, it can be shown that $\psi_{2,1}(t)$ (and hence $\psi_{2,1}^{*}(t)$) is non-decreasing (non-increasing) in $t$. Now to show that $\lim_{t \to \infty} \; \psi_{2,1}^{*}(t)=c_{0,1}$, we will show that $\psi_{2,1}(t)\geq\; (\leq)\; \psi_{0,1}(t),\; \forall \;t$. Let us fix $t$, then 
$  l_1(\psi_{0,1}(t)|t)=0 $	and	$l_2(\psi_{2,1}(t)|t)=0.$
\vspace*{2mm}

Let $s_0=\frac{1}{\psi_{0,1}(t)},\; M(s)=s\,W^{'}(\frac{s}{s_0})f_1(s),\; M_1(s)=h_1(t|s)$ and $M_2(s)=H_1(t|s),\; s\in \Re_{++}.$ Using hypotheses of the corollary, Lemma 3.1 and Lemma 2.1, we conclude that
\begin{align*}
	\qquad \quad \; 	l_2(\psi_{0,1}(t)|t) &\leq\; (\geq)\; \frac{h_1(t|s_0)}{H_1(t|s_0)}\; 	l_1(\psi_{0,1}(t)|t)
	=0
\end{align*}
Since $l_2(c|t)$ is an non-decreasing function of c (using $A_1$) and $\psi_{2,1}(t)$ is the unique solution of $l_2(c|t)=0$, we conclude that $\psi_{0,1}(t)\leq \;(\geq) \; \psi_{2,1}(t)$. Hence $c_{0,1}=\lim_{t \to \infty} \psi_{0,1}(t)$ $\leq\; (\geq) \; \lim_{t \to \infty}\psi_{2,1}(t)$ and $\lim_{t \to \infty} \psi_{2,1}^{*}(t)= \min\{c_{0,1},\lim_{t \to \infty} \psi_{2,1}(t)\}$ $=c_{0,1}$ ($\lim_{t \to \infty} \psi_{2,1}^{*}(t)\\= \max\{c_{0,1},\lim_{t \to \infty}\psi_{2,1}(t)\}=c_{0,1}$). Note that $\psi_{2,1}^{*} (t)\geq \; (\leq)\; \psi_{0,1}(t),\; \forall\;t$. Since $l_1(c|t)$ is a non-decreasing function of $c$, we have 
$$l_1(\psi_{2,1}^{*}(t)|t)\geq\;(\leq)\; l_1(\psi_{0,1}(t)|t)=0,\; \forall\;t.$$
Hence the result follows.
\\~\\ The proof of part (ii) is immediate using Theorem 3.1.1 and the fact that $l_1(c|t)$ is a non-decreasing function of $c$.
\end{proof}

It is straightforward to see that the estimator $\delta_{\psi_{0,1}},$ defined in Corollary 3.1.1 (i), is the generalized Bayes estimator with respect to the non-informative density $\pi(\theta_1,\theta_2)=\frac{1}{\theta_1 \theta_2},\;(\theta_1,\theta_2)\in\Theta_0$.

\subsection{\textbf{Improvements Over the BSEE of $\theta_2$}}
\label{sec:3.2}
\setcounter{equation}{0}
\renewcommand{\theequation}{3.2.\arabic{equation}}

\noindent
\vspace*{2mm}

As proofs of various results stated in this section are similar to the proofs of similar results of the last section, they are being omitted.
The following theorem provides a class of estimators that improve upon the BSEE $\delta_{c_{0,2},2}(\underline{X})=c_{0,2} X_2$, where $c_{0,2}$ is the unique solution of the equation
$\int_{-\infty}^{\infty} z\,W^{'}(cz)f_2(z) \,dz=0.$
\\~\\\textbf{Theorem 3.2.1.} Suppose that, for any fixed $\Delta\geq 1$ and $t$, $H_2(\frac{t}{\Delta}\vert s)/H_2(t\vert s)$ is non-increasing (non-decreasing) in $s\in S_2$. Consider a scale equivariant estimator $\delta_{\psi_2}(\underline{X})=\psi_2(D)X_2$ for estimating $\theta_2$, where $\lim_{t\to\infty} \psi_2(t)=c_{0,2}$, $\psi_2(t)$ is an non-increasing (non-decreasing) function of $t$ and $\int_{0}^{\infty}s\, W^{'}(\psi_2(t)s)\; H_2(t\vert s)\,f_2(s)ds \, \leq \,(\geq) \,0,\; \forall\;t$. Then, $\forall\; \underline{\theta}\in \Theta_0$, the estimator $\delta_{\psi_2}(\underline{X})$ dominates the BSEE $\delta_{c_{0,2},2}(\underline{X})=c_{0,2}X_2$.\vspace*{2mm}
\\~\\ The following corollary provides the B-Z type improvements over the BSEE $\delta_{c_{0,2},2}(\underline{X})=c_{0,2} X_2$.
\\~\\	\textbf{Corollary 3.2.1. (i)} Suppose that assumptions of Theorem 3.1.1 hold. Further suppose that, for every fixed $t$, the equation
$$l_3(c|t)=\int_{0}^{\infty} \;s\, W^{'}(cs)\; H_2(t|s) f_2(s)\;ds =0$$
has the unique solution $c\equiv \psi_{0,2}(t)$. Then, the estimator $\delta_{\psi_{0,2}}(\underline{X})=\psi_{0,2}(D)X_2$ improves upon the BSEE $\delta_{c_{0,2},2}(\underline{X})=c_{0,2}X_2$, $\forall\; \underline{\theta}\in\Theta_0$.
\\~\\\textbf{(ii)} In addition to assumptions of (i) above, suppose that $\psi_{1,2}:\Re_{++} \rightarrow \Re $ is such that $\psi_{1,2}(t) \leq \, (\geq) \; \psi_{0,2}(t), \;  \forall \;t, \; \psi_{1,2}(t) $ is non-increasing (non-decreasing) in $t$ and $\lim_{t \to \infty}\; \psi_{1,2}(t)\\ = c_{0,2}$.
Then, $\forall\; \underline{\theta}\in\Theta_0$, the estimator $\delta_{\psi_{1,2}}(\underline{X})=\psi_{1,2}(D)X_2$ dominates the BSEE $\delta_{c_{0,2},2}(\underline{X})=c_{0,2}X_2$.
\vspace{2mm}

In the following corollary we provide the Stein type improvements over the BSEE $\delta_{c_{0,2},2}(\underline{X}).$
\\~\\ \textbf{Corollary 3.2.2. (i)} Suppose that for any fixed $\Delta\geq 1$ and $t$, $h_2\left(\frac{t}{\Delta}\vert s\right)/h_2(t\vert s)$ is non-increasing (non-decreasing) in $s\in S_2$ and let $\psi_{0,2}(t)$ be as defined in Corollary 3.2.1. In addition suppose that, for any $t$, the equation
$$l_4(c|t)=\int_{0}^{\infty} \;s\, W^{'}(c\,s)\; h_2(t|s) f_2(s)\;ds =0$$
has the unique solution $c\equiv \psi_{2,2}(t)$. Let $\psi_{2,2}^{*}(t)= \max\{c_{0,2},\psi_{2,2}(t)\}$ ($\psi_{2,2}^{*}(t)= \min\{c_{0,2},\psi_{2,2}(t)\}$) and $\delta_{\psi_{2,2}^{*}}(\underline{X}) =  \psi_{2,2}^{*}(D)X_2$. Then
$$R_2(\underline{\theta},\delta_{\psi_{2,2}^{*}})\leq R_2(\underline{\theta},\delta_{c_{0,2},2})\, \;\;\; \forall \; \; \underline{\theta} \in \Theta_0.$$
\textbf{(ii)} In addition to assumptions of (i) above, suppose that $\psi_{3,2}:\Re_{++} \rightarrow \Re $ is such that $\psi_{3,2}(t) \leq \; (\geq) \; \psi_{2,2}(t), \; \; \forall \;t$ and $\psi_{3,2}(t) $ is non-increasing (non-decreasing) in $t$. For fixed $t$, define $\psi_{3,2}^{*}(t)=\max\{c_{0,2},\psi_{3,2}(t)\}\;(\psi_{3,2}^{*}(t)=\min\{c_{0,2},\psi_{3,2}(t)\})$ and $\delta_{\psi^{*}_{3,2}}(\underline{X})=\psi^{*}_{3,2}(D)X_2$.
Then 
$$R_2(\underline{\theta}, \delta_{\psi^{*}_{3,2}}) \leq R_2(\underline{\theta}, \delta_{c_{0,2},2}),\; \; \forall \; \underline{\theta} \in \Theta_0.$$

It is easy to verify that the B-Z type estimator $\delta_{\psi_{0,2}}(\cdot)$, derived in Corollary 3.2.1 (i), is the generalized Bayes estimator with respect to the non-informative prior density $\pi(\theta_1,\theta_2)=\frac{1}{\theta_1 \theta_2},\;\;\;\;(\theta_1,\theta_2)\in\Theta_0.$\vspace*{2mm}

The results of Theorems 3.1.1-3.2.1 (or Corollaries 3.1.1-3.1.2 and Corollaries 3.2.1-3.2.2) are applicable to various studies carried out in the literature for specific bivariate probability models, having independent marginals, and specific loss function (e.g., Misra and Dhariyal (\citeyear{MR1326266}), Vijayasree et al. (\citeyear{MR1345425}), etc.). These results also extend the study of Kubokawa and Saleh (\citeyear{MR1370413}) to general bivariate scale models.
\\~\\ Now we provide an application of the results derived in subsections 3.1-3.2 to a situation where results of Kubokawa and Saleh (\citeyear{MR1370413}) are not applicable.

\subsection{\textbf{Applications}}
\label{sec:3.3}

\setcounter{equation}{0}
\renewcommand{\theequation}{3.3.\arabic{equation}}

\noindent
\vspace*{2mm}

In the following example, we consider a bivariate model due to Cheriyan and Ramabhadran's (see Kotz et al. \citeyear{MR1788152}) and study estimation of order restricted scale parameters.
\\~\\ \textbf{Example 3.3.1.}  Let $X_1$ and $X_2$ be two dependent random variables with joint pdf \eqref{eq:1.2},
where $\underline{\theta}=(\theta_1,\theta_2)\in\Theta_0$ and \vspace*{1.5mm}

\begin{center}$ f(z_1,z_2)= \begin{cases} e^{-z_1} (1-e^{-z_2}),\; \; 0<z_2<z_1\\
	e^{-z_2} (1-e^{-z_1}),\; \; 0<z_1<z_2\\
	0,\qquad \qquad \qquad  \text{    otherwise}
\end{cases}.   $\end{center}

The above bivariate distribution is a special case of Cheriyan and Ramabhadran's bivariate gamma distribution (see Kotz et al. \citeyear{MR1788152}). Here random variable $X_i$ follows Gamma distribution with pdf $f_i(\frac{x}{\theta_i})=\frac{x}{\theta_i^2}\,e^{-\frac{x}{\theta_i}},\; x>0$, $i=1,2$.
\\~\\For estimation of $\theta_i,\;i=1,2$, consider the squared error loss function
$ L_i(\underline{\theta},a)=\left(\frac{a}{\theta_i}-1\right)^2,$ $ \underline{\theta}\in \Theta_0,\;a\in \Re_{++},\;i=1,2.$
The BSEE of $\theta_i$ is $\delta_{c_{0,i},i}(\underline{X})=\frac{1}{3} X_i$, $i=1,2$ ($c_{0,i}=\frac{1}{3},\;i=1,2$). We have $S_1=S_2=[0,\infty)$.
\\~\\\textbf{\underline{Estimation of $\theta_1$}:} \vspace*{2mm}

\noindent
For any $s\in S_1$, the pdf and df of $Z_{s}^{(1)}$, respectively, are
\begin{align*}
h_1(t|s)=\begin{cases} 1-e^{-st} ,& 0<t< 1\\ e^{-s(t-1)}(1-e^{-s}),& 1\leq t<\infty \\0,&\text{otherwise} \end{cases}
\;\;\text{and } \;
H_1(t|s)=\begin{cases}0,& t<0\\ t-\frac{1}{s}+\frac{e^{-st}}{s} ,&0\leq t< 1\\ \frac{(1-e^{-s})(1-e^{-s(t-1)})}{s},& t\geq 1 \end{cases}.
\end{align*}
It is easy to see that, for any fixed $\Delta\geq 1$ and $t$, $h_1\!\left(\frac{t}{\Delta}\vert s\right)/h_1(t\vert s)$ (and hence $H_1\!\left(\frac{t}{\Delta}\vert s\right)/H_1(t\vert s)$) is non-decreasing in $s\in S_1=\Re_{++}$. We have
\begin{align*}
\psi_{0,1}(t)&= \frac{\int_{-\infty}^{\infty} s\,H_1(t|s)f_1(s)ds}{\int_{-\infty}^{\infty}s^2 H_1(t|s)f_1(s)ds}
= \begin{cases} \frac{2t-1+\frac{1}{(t+1)^2}}{6t-2+\frac{2}{(t+1)^3}}, \qquad 0<t<1 \\~\\ \frac{2-\frac{1}{t^2}+\frac{1}{(t+1)^2}}{6-\frac{2}{t^3}+\frac{2}{(t+1)^3}},\qquad t\geq 1 \end{cases},\\	
\psi_{2,1}(t)&= \frac{\int_{-\infty}^{\infty} s\,h_1(t|s)f_1(s)ds}{\int_{-\infty}^{\infty}s^2 h_1(t|s)f_1(s)ds}
= \begin{cases} \frac{1}{3}\,\frac{1-\frac{1}{(t+1)^3}}{1-\frac{1}{(t+1)^4}}, \qquad \;\; 0<t<1 \\~\\ \frac{1}{3}\frac{\frac{1}{t^3}-\frac{1}{(t+1)^3}}{\frac{1}{t^4}-\frac{1}{(t+1)^4}},\qquad \;\;\;t\geq 1 \end{cases}\\
\text{and}\qquad
\psi_{2,1}^{*}(t)&=\min\{c_{0,1}\,,\psi_{2,1}(t)\}=\begin{cases} \frac{1}{3}\,\frac{1-\frac{1}{(t+1)^3}}{1-\frac{1}{(t+1)^4}}, \qquad \qquad\qquad\;\; 0<t<1 \\~\\
	\frac{1}{3} \min\Bigg\{1,\frac{\frac{1}{t^3}-\frac{1}{(t+1)^3}}{\frac{1}{t^4}-\frac{1}{(t+1)^4}}  \Bigg\},\qquad\; t\geq 1
\end{cases}.
\end{align*}
Here $\psi_{0,1}(t)$ and $\psi_{2,1}(t)$ are non-decreasing in $t\in \Re_{++}$ and $\lim_{t \to \infty} \psi_{0,1}(t)=\frac{1}{3}=c_{0,1}$.

\noindent
Using Corollary 3.1.1 (i), the B-Z type estimator dominating the BSEE $\delta_{c_{0,1},1}(\underline{X})=\frac{1}{3} X_1$ is 
\begin{equation}\label{eq:3.2.1}
\delta_{\psi_{0,1}}(\underline{X})=\psi_{0,1}(D)X_1= \begin{cases} \frac{2D-1+\frac{1}{(D+1)^2}}{6D-2+\frac{2}{(D+1)^3}}X_1, \qquad X_1>X_2 \\~\\ \frac{2-\frac{1}{D^2}+\frac{1}{(D+1)^2}}{6-\frac{2}{D^3}+\frac{2}{(D+1)^3}}X_1,\qquad X_1\leq X_2 \end{cases}.	
\end{equation}
Here $\delta_{\psi_{0,1}}(\cdot)$ is also the generalized Bayes estimator with respect to the non-informative prior density on $\Theta_0$. 

\noindent
Using Corollary 3.1.2 (i), the Stein type estimator dominating the BSEE $\delta_{c_{0,1},1}(\underline{X})=\frac{1}{3} X_1$ is 
\begin{equation}\label{eq:3.2.2}
\delta_{\psi_{2,1}^{*}}(\underline{X})=\psi_{2,1}^{*}(D)X_1=\begin{cases} \frac{1}{3}\frac{1-\frac{1}{(D+1)^3}}{1-\frac{1}{(D+1)^4}}X_1\,,  \qquad \qquad\qquad  X_1> X_2 \\ \frac{1}{3}  \min\bigg\{1,\frac{\frac{1}{D^3}-\frac{1}{(D+1)^3}}{\frac{1}{D^4}-\frac{1}{(D+1)^4}} \bigg\}X_1  \, ,\quad   X_1\leq X_2 \end{cases}.	
\end{equation}
\textbf{\underline{Estimation of $\theta_2$}:} \vspace*{2mm} 

\noindent
For any $s\in S_2=[0,\infty)$, the pdf and df of $Z_{s}^{(2)}$, respectively, are
\begin{align*}
h_2(t|s)&=\begin{cases} \frac{e^{-\frac{s}{t}}\,e^{s}(1-e^{-s})}{t^2} ,\; \; \; \,0<t< 1\\ \frac{(1-e^{-\frac{s}{t}})}{t^2},\quad\qquad 1\leq t<\infty\\0,\qquad\qquad\qquad\; \text{elsewhere} \end{cases}
\text{and } \;
H_2(t|s)&=\begin{cases}0,& t\leq 0\\ \frac{e^{-s\left(\frac{1}{t}-1\right)}-e^{-\frac{s}{t}}}{s} ,&0<t< 1\\ 1-\frac{1}{t}+\frac{1}{s}-\frac{e^{-\frac{s}{t}}}{s},& t\geq 1 \end{cases}.
\end{align*}
One can easily see that, for any fixed $\Delta\geq 1$ and $t$, $h_2\!\left(\frac{t}{\Delta}\vert s\right)/h_2(t\vert s)$ (and hence $H_2\!\left(\frac{t}{\Delta}\vert s\right)/H_2(t\vert s)$) is non-increasing in $s\in (0,\infty)$.
Let $\psi_{0,2}(\cdot)$ and $\psi_{2,2}(\cdot)$ be as defined in Corollaries 3.2.1 (i) and 3.2.2 (i), respectively, so that for fixed $t$, we have
\begin{align*}
\psi_{0,2}(t)&= \frac{\int_{-\infty}^{\infty} s\,H_2(t|s)f_2(s)ds}{\int_{-\infty}^{\infty}s^2 H_2(t|s)f_2(s)ds}
= \begin{cases} \frac{t}{2}\frac{1-\frac{1}{\left(1+t\right)^2}}{1-\frac{1}{\left(1+t\right)^3}}, & 0<t<1 \\~\\ \frac{3-\frac{2}{t}-\frac{t^2}{\left(1+t\right)^2}}{8-\frac{6}{t}-\frac{2t^3}{\left(1+t\right)^3}},& t\geq 1 \end{cases},\\	
\psi_{2,2}(t)&= \frac{\int_{-\infty}^{\infty} s\,h_2(t|s)f_2(s)ds}{\int_{-\infty}^{\infty}s^2 h_2(t|s)f_2(s)ds}
= \begin{cases} \frac{1}{3t}\frac{1-\frac{1}{\left(1+t\right)^3}}{1-\frac{1}{\left(1+t\right)^4}}, & 0<t<1 \\~\\ \frac{1}{3}\frac{1-\frac{t^3}{\left(1+t\right)^3}}{1-\frac{t^4}{\left(1+t\right)^4}},& t\geq 1 \end{cases}\\
\text{and}\qquad
\psi_{2,2}^{*}(t)&=\max\{c_{0,2}\,,\psi_{2,2}(t)\}=\begin{cases} \frac{1}{3}\,\max\bigg\{1, \frac{1}{t}\frac{1-\frac{1}{\left(1+t\right)^3}}{1-\frac{1}{\left(1+t\right)^4}}\bigg \}, & 0<t<1 \\~\\
	\frac{1}{3},& t\geq 1
\end{cases}.
\end{align*}
Here $\psi_{0,2}(t)$ and $\psi_{2,2}(t)$ are non-increasing in $t\in \Re_{++}$ and $\lim_{t \to \infty} \psi_{0,2}(t)=\frac{1}{3}$.

\noindent
Using Corollary 3.2.1 (i), the B-Z type estimator dominating the BSEE $\delta_{c_{0,2},2}(\underline{X})=\frac{1}{3} X_2$ is 
$$\delta_{\psi_{0,2}}(\underline{X})=\psi_{0,2}(D)X_2= \begin{cases} \frac{D}{2}\frac{1-\frac{1}{\left(1+D\right)^2}}{1-\frac{1}{\left(1+D\right)^3}}X_2, & X_1>X_2 \\~\\ \frac{3-\frac{2}{D}-\frac{D^2}{\left(1+D\right)^2}}{8-\frac{6}{D}-\frac{2D^3}{\left(1+D\right)^3}}X_2,& X_1\leq X_2 \end{cases}	. $$
Using Corollary 3.2.2 (i), the Stein type estimator dominating the BSEE $\delta_{c_{0,2},2}(\underline{X})=\frac{1}{3} X_2$ is 
$$\delta_{\psi_{2,2}^{*}}(\underline{X})=\psi_{2,2}^{*}(D)X_2=\begin{cases} \frac{1}{3}\,\max\bigg\{1, \frac{1}{D}\frac{1-\frac{1}{\left(1+D\right)^3}}{1-\frac{1}{\left(1+D\right)^4}}\bigg \}X_2\,,  & X_1> X_2 \\  \frac{1}{3}X_2\, ,& X_1\leq X_2 \end{cases}.$$

\subsection{\textbf{Simulation Study For Estimation of Scale Parameter $\theta_1$}}
\label{sec:3.4}
\noindent
\vspace*{2mm}

In Example 3.2.1, we have considered a Cheriyan and Ramabhadran's bivariate gamma distribution with unknown order restricted scale parameters (i.e., $\theta_1\leq \theta_2$). To further evaluate the performances of various estimators of $\theta_1$ under the squared error loss function, in this section, we compare the risk performances of the  BSEE $\delta_{c_{0,1},1}(\underline{X})=\frac{X_1}{3}$, the B-Z estimator $\delta_{\psi_{0,1}}$ and the Stein (1964) type estimator $\delta_{\psi^{*}_{2,1}}$ (as defined in \eqref{eq:3.2.1} and \eqref{eq:3.2.2}), numerically, through Monte Carlo simulations. The simulated risks of the BSEE, the B-Z estimator and the Stein estimator have been computed based on 10000 simulations from relevant distributions.

\begin{figure}[h!]
\centering
% include first image
\includegraphics[width=140mm,scale=1.5]{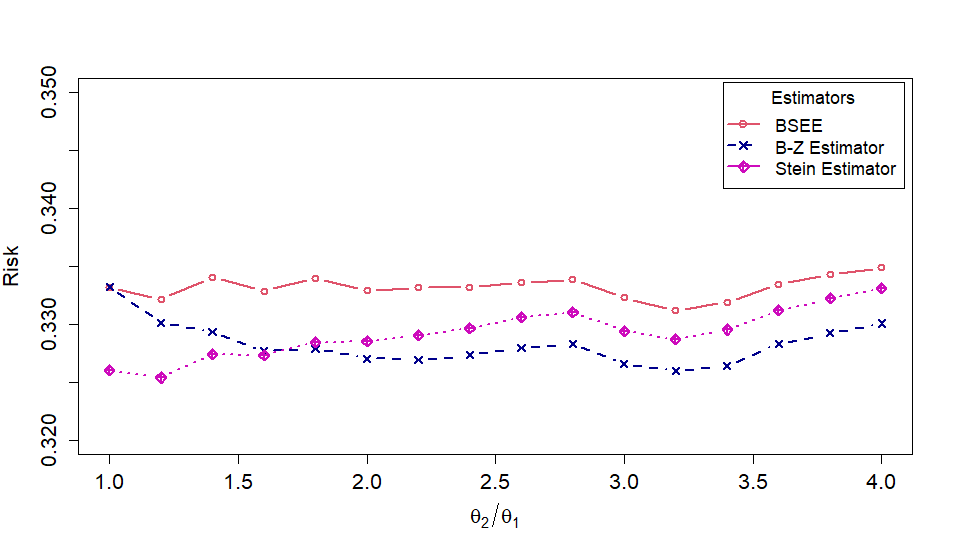} 
\caption{Risk plot of $\delta_{c_{0,1},1}$ (BSEE), $\delta_{\psi_{0,1}}$ (B-Z type estimator) and $\delta_{\psi^{*}_{2,1}}$ (Stein type estimator) estimators against the values of $\frac{\theta_2}{\theta_1}$.}
\label{fig8}
\end{figure}

The simulated values of risks of various estimators are plotted in Figure \ref{fig8}. The following observations are evident from Figure \ref{fig8}:
\\(i) The B-Z type and the Stein type estimators always perform better than the BSEE, which 
\FloatBarrier

\noindent
is in conformity with theoretical findings of Example 3.2.1.
\\(ii) There is no clear cut winner between the B-Z type estimator $\delta_{\psi_{0,1}}$ and the Stain type estimator $\delta_{\psi^{*}_{2,1}}$. The Stein type estimator performs better than the B-Z type estimator, for small values of $\frac{\theta_2}{\theta_1}$, and the B-Z type estimator dominates the Stein type estimator for the large values of $\frac{\theta_2}{\theta_1}$. 

\section{\textbf{Concluding Remarks}}

The problem of estimation of order restricted location/scale parameters is widely studied for specific probability models, having independent marginals, and specific loss functions. In this paper, we unify these studies by considering a general bivariate location/scale model and a general loss function. We drive a class of estimators dominating over BLEE/BSEE using the IERD approach of Kubokawa (\citeyear{MR1272084}). We also obtain the Brewster-Zidek (\citeyear{MR381098}) type and the Stein (\citeyear{MR171344}) type estimators that dominate the BLEE/BSEE under the general loss function. We also demonstrate applications of our results to two bivariate probability models which have not been studied in the literature.

\section*{\textbf{Funding}}

This work was supported by the [Council of Scientific and Industrial Research (CSIR)] under Grant [number 09/092(0986)/2018].

\bibliographystyle{apalike}
\bibliography{Paper1}

\end{document}